\documentclass[11pt]{article}
	
	\newcommand{\blind}{0}
    \makeatletter
    \renewcommand\section{\@startsection {section}{1}{\z@}%
                                       {-3.5ex \@plus -1ex \@minus -.2ex}%
                                       {2.3ex \@plus.2ex}%
                                       {\normalfont\fontfamily{phv}\fontsize{16}{19}\bfseries}}
    \renewcommand\subsection{\@startsection{subsection}{2}{\z@}%
                                         {-3.25ex\@plus -1ex \@minus -.2ex}%
                                         {1.5ex \@plus .2ex}%
                                         {\normalfont\fontfamily{phv}\fontsize{14}{17}\bfseries}}
    \renewcommand\subsubsection{\@startsection{subsubsection}{3}{\z@}%
                                        {-3.25ex\@plus -1ex \@minus -.2ex}%
                                         {1.5ex \@plus .2ex}%
                                         {\normalfont\normalsize\fontfamily{phv}\fontsize{14}{17}\selectfont}}
    \makeatother
	\usepackage{amsmath}
         \usepackage{amsfonts}
	\usepackage{graphicx}
	\usepackage{enumerate}
	\usepackage{xcolor}
    \usepackage{makecell}  
	\usepackage{natbib} %comment out if you do not have the package
	\usepackage{url} % not crucial - just used below for the URL
	\def\ZZ{\hbox{\sf Z\kern-.43emZ}}
\def\RR{\hbox{\sf I\kern-.14em\hbox{R}}}%%%%%%%%%%%%%%%%%%%%%%%%%%%%%%%%%%%%%%%%%%%%%%%%%%%%%%%%%%%%%%%%%%%%%%%
	
    \usepackage{color}
    \definecolor{EditBlue}{RGB}{0,128,255}
    
    \usepackage{endnotes}
    \let\footnote=\endnote
    \let\enotesize=\normalsize
    \def\notesname{Endnotes}%
    \def\makeenmark{$^{\theenmark}$}
    \def\enoteformat{\rightskip0pt\leftskip0pt\parindent=1.75em
      \leavevmode\llap{\theenmark.\enskip}}
    \usepackage{algorithm}
    \usepackage[noend]{algorithmic}
    \usepackage{enumitem}
    \usepackage{booktabs}
    \usepackage{setspace}
    \usepackage[font=singlespacing]{caption}
    \usepackage{rotating}
    \usepackage{caption}

    \DeclareCaptionLabelFormat{andtable}{#1~#2  \&  \tablename~\thetable}

    \usepackage{booktabs}
    \usepackage{multirow}
    \usepackage{tabularx,array}
    \newcolumntype{Y}{>{\centering\arraybackslash}X}    

\begin{document}
		
			%%%%%%%%%%%%%%%%%%%%%%%%%%%%%%%%%%%%%%%%%%%%%%%%%%%%%%%%%%%%%%%%%%%%%%%%%%%%%%
		\def\spacingset#1{\renewcommand{\baselinestretch}%
			{#1}\small\normalsize} \spacingset{1}
		%%%%%%%%%%%%%%%%%%%%%%%%%%%%%%%%%%%%%%%%%%%%%%%%%%%%%%%%%%%%%%%%%%%%%%%%%%%%%%
		
		\if0\blind
		{
			\title{\bf {On the Optimization of Benefit to Cost Ratios for \\ Public Sector Decision Making}}
			\author{Frederick ``Forrest'' Miller $^a$ , Yaren Bilge Kaya $^b$ , Geri Louise Dimas $^c$, \\ Renata A. Konrad $^d$, Kayse Lee Maass $^a$, Andrew C. Trapp $^{f}$\\ \\
			$^a$ Department of Mechanical and Industrial Engineering, Northeastern University, \\ Boston, MA, USA \\
             $^b$ Department of Industrial Engineering and Operations Research, \\ Columbia University, New York City, NY, USA \\
             $^c$ Department of Information Systems and Analytics, \\ Bryant University, Smithfield, RI, USA\\
             $^d$  Department of Industrial and Systems Engineering, \\ Rensselaer Polytechnic Institute, Troy, NY, USA \\
% $^e$ Department of Mechanical and Industrial Engineering, \\ Northeastern University, Boston, MA, USA \\
$^f$ Data Science Program, Worcester Polytechnic Institute, Worcester, MA, USA}

			\date{}
			\maketitle
		} \fi
		
		\if1\blind
		{

            \title{\bf \emph{IISE Transactions} \LaTeX \ Template}
			\author{Author information is purposely removed for double-blind review}
			
\bigskip
			\bigskip
			\bigskip
			\begin{center}
				{\LARGE\bf Optimizing Benefit-to-Cost Ratios: \\ A Framework for Public Sector Decisions}
			\end{center}
			\medskip
		} \fi
		% \bigskip
            \vspace{-3mm}
  
	\begin{abstract}
Decision making in the public sector emphasizes objectives like equity and efficiency alongside standard economic returns to reflect the broader cares of society and public beneficiaries.
These contexts often feature multiple stakeholders and objectives in complex environments, requiring the balancing of operations, deployment and change management. 
Cost-benefit analysis is a prevailing decision-making paradigm in the public sector, employing tools like the benefit to cost ratio (BCR) to evaluate alternatives (actions) beyond the status quo.
Yet, ex ante identification of alternatives can be challenging and may miss opportunities.
We introduce a new framework for public sector decision making that maximizes the BCR. 
Our framework evaluates numerous actions, subject to original conditions and \textit{action constraints} that ensure movement beyond the status quo. 
We represent this problem as a mixed-integer linear fractional program (MILFP) and employ Newton's Method to convert the MILFP to a sequence of linearized mixed-integer programs, retaining tractability for large instances. 
We demonstrate our framework on a New York City (NYC) runaway and homeless youth shelter system case study. 
Our optimization-based framework yields data-informed recommendations for NYC shelter and service expansion decisions. 
Our framework enables iteration and comparison across action constraints of increasing degree, empowering effective assessment of marginal activity.

	\end{abstract}
			
	\noindent%

	%\newpage
	\spacingset{1.5} % DON'T change the spacing!

\section{Introduction} \label{intro}
\looseness-1
%Like many decision-making contexts, p
Investments in public sector activities are influenced by benefits and costs in ways that differ from strict profit maximization, instead valuing benefit to cost \textit{ratios} \citep{zerbe2012primer,whittington1984economic}.
Standard accounting techniques for obtaining public sector costs exist, but accurate estimations of benefits are often challenging and nuanced \citep{Gass-1994}.
Benefits can accrue to a variety of stakeholders in ways that may be less obvious than fixed and operational costing estimates.
Decision-theoretic results support the superiority of benefit to cost ratios (BCRs) in evaluating gains and losses to make informed decisions where there is aversion to risk, so familiar in the public sector~\citep{bangforbuck}.
Indeed, in settings with significant resource scarcity for investing in critical societal needs, simply minimizing costs, or maximizing the difference between benefits and costs while consuming no more than a budget, may be shortsighted.

Organizations that operate in complex and dynamic contexts are often sensitive to instability in their operations.
These contexts are influenced by a variety of factors including, pressing demands, multiple stakeholder interests, and regulatory expectations.
Public sector decisions, in particular, have far-reaching consequences that impact the lives of citizens and communities, and the effects of abstractly modeling and optimizing such decisions  must be balanced with the pragmatic reality.
If an optimal solution is identified, it may be impractical to deploy due to the associated change management.
In such contexts, a road map that charts out the next best steps is needed for organizational decision-making.

One public sector context featuring heightened resource limitations is the provision of critical support services for vulnerable individuals.
%Vulnerable populations have needs that, in the collective, tend to greatly outstrip supply.
%The organizations that provide services to these groups are vital to ensure they avoid being overlooked by society.
In particular, the homeless, including runaway and homeless youth (RHY), are served by a variety of public sector organizations that provide programming, shelters, and housing.
These organizations offer a safe place to stay and myriad services, including education, health care, employment, and vocational training to help RHY become self-sufficient and find permanent housing. 
Effective prevention and support services have been shown to reduce vulnerability and recidivism for unstably housed youth \citep{ltmurphy,clawson, survivingStreetsNyc}.
% In New York City, it is estimated that annually \textbf{5,734} youth are on the street, while collective youth-specific shelter bed capacity runs around \textbf{1,195} \citep{Morton-2019}.
While RHY with unaddressed needs are at high risk of exploitation, including human trafficking (HT) \citep{Lwilliams,Middleton}, the supply of existing beds in the United States (U.S.) shelter system falls woefully short of the collective demand for RHY seeking shelter services~\citep{structurePractice}.

%Human trafficking (HT) is a form of modern slavery and can be described as a series of criminal events exploiting humans for benefits \citep{doi:10.1177/1363461507081635}.
%The 2021 global estimate indicates that nearly 50 million people experience some form of modern slavery on any given day \citep{ILOModernSlaveryEstimates_2021} and annually generate more than 150 billion USD in illicit gains globally \citep{htImpactEstimate}.
%Traffickers are known to target individuals with vulnerabilities such as those that contribute to youth homelessness \citep{Curtis-2008}.

%PARAGRAPH 6: CONTRIBUTIONS
\looseness-1
We propose a new framework for evaluating a potentially vast array of decision alternatives beyond maintaining the status quo in the public sector domain.
Our framework identifies a (societal) benefit to (public sector investment) cost ratio maximizing the next best course of action, naturally seeking alternatives that yield the largest benefit per marginal activity from the status quo.
To the best of our knowledge, we are the first to put forth a decision-making framework applicable to a broad set of contexts where many alternatives from the status quo must be evaluated on the basis of maximizing the BCR, and specifically for public sector decision making. 

In contrast to evaluating a small, finite number of user-defined alternatives, we propose a new algorithmic framework that implicitly represents evaluating combinatorially many alternatives to identify one that maximizes the BCR objective.

As a case study, we investigate resource allocation decisions for programs to support vulnerable RHY that are at high risk of HT within New York City (NYC). Specifically, we provide guidance on how to expand capacity within NYC to meet the shelter and supportive service needs of RHY in a ratio-maximizing manner. We estimate benefits accrued for RHY through shelter and service provision versus the costs of expanding resources in existing organizations, building new shelters, and assigning youth to those shelters. The problem of determining where to locate and open shelters, as well as how to deploy additional service capacity, differs from that of only determining how to place youth within available shelters~\citep{OriginalModel}; our approach considers both shelter expansion, as well as youth placement and service allocation. 
%By considering marginal activity from the status quo, our approach finds the cost-minimizing allocation of resources satisfying all expected RHY demand.
%In particular, we allow for actions (e.g., determine which shelters to open, where to open them, and what their capacities should be) involving decisions for a variety of new shelters that fit the RHY profile needs informed by the data-driven capacity needs observed from \cite{OriginalModel}.

%This work addresses a problem with high societal need and contributes to informing access and social welfare of RHY in NYC.

%Through this we investigate optimizing the BCR versus a set of potential competing objective functions, and demonstrate the advantages of the ratio optimization in best utilizing expanded resources as compared to alternative objectives.

The remainder of this paper is organized as follows. Section \ref{background} covers relevant background. Section \ref{Methods} introduces our new algorithmic framework. Section \ref{Model} introduces our case study on expanding RHY shelter system capacity in NYC, detailing the mathematical models, algorithms, and metrics developed. Section \ref{CompuationalExperiments} discusses our computational experiments and findings concerning optimizing the BCR objective for public sector capacity expansion. Section \ref{keyTakeaways} concludes.

\section{Background}\label{background}
% OR \ rhy & homeless youth
In this section, we discuss relevant background for our study, including current practices for public sector decision making and related applications of fractional programming.
\subsection{\emph{Decision Making in the Public Sector}}\label{DMPS}
%PARAGRAPH 1
\looseness-1
% The public sector centers on delivering resources and services for the common good.
% This contrasts with the private sector that focuses on returning value for shareholders.
% The public sector also differs from the nonprofit sector, which aims to address social, environmental, cultural, and educational needs of communities that are underserved by society \citep{IndianaNonprofitsProject_WhatIsANonprofit}. 
% It also is distinct from fourth sector movements such as social enterprises that play a hybrid role of addressing societal needs while pursuing financial returns.
% Within the public sector jurisdiction lie urban services and transportation; the natural environment, natural resources and land conservation; public utilities; national defense; education; as well as crime and justice.
% %Crime: Blumstein, A. (2007) "An OR Missionary's Visits to the Criminal Justice System," Ops. Res., 55, 1: 14-23.
% %Narcotics: Caulkins, J. P. (2005) How Goes the War on Drugs? An Assessment of US Drug Problems and Policy,
% % Decisions made in the public sector are more likely than those in the private sector to impact those beyond actual decision makers.
% While financials are one important criterion in the public sector, they are joined by a broader set of objectives and goals that characterize society and the beneficiaries of public services -- equity and efficiency \citep{Gass-1994}.
%==============

%PARAGRAPH 2
Public sector decision making (PSDM) typically involves multiple stakeholders with conflicting objectives. Some stakeholders have objectives that are often less obvious and quantifiable than their financial counterparts, yet are critical in evaluating the merit of a proposed intervention.
Informed decisions, particularly in PSDM, necessitate understanding and aligning disparate views and values to arrive at tractable problem representations.
Community-based operations research \citep{johnson2007community,johnson2012community}, soft operations research \citep{mingers2011soft}, and problem structuring methods \citep{ackermann2012problem,smith2019characteristics} promote process-oriented means to capture the values and perceptions of multiple stakeholders regarding their goals and alternative identification, facilitating improved model development, and ultimately implementation.
Similarly, multiple-criteria decision making \citep{greco2016multiple} enables decision makers to evaluate alternatives with respect to diverse objectives.
The knowledge of performance in multiple public sector dimensions can lead to more informed decisions.
%==============

%PARAGRAPH 3

A predominant decision-making framework in PSDM for evaluating a public good intervention beyond the \textit{status quo} is cost-benefit analysis (CBA), and similarly, cost-effectiveness analysis (CEA) in the military, environmental and health economics contexts \citep{boardman2017cost}.
In CBA, the \textit{benefit to cost} (also \textit{cost-benefit}) ratio (BCR) measures the present value of benefits the public accrues to that of governmental cost outlay.
The BCR is an indicator that evaluates whether the likely benefits to expenditures generated are superior to the next best alternative or status quo.
In complex decision contexts, such as those often found in the public sector, incremental improvements may be favored over drastic changes that represent a significant departure from the status quo.
In the mid-twentieth century, BCRs began to be used in water allocation and land management contexts in the U.S. and are now widely used in governmental agencies worldwide to evaluate resource allocation decisions \citep{fuguitt1999cost}.

\looseness-1
The BCR is related to the notion of marginal utility, which captures the change in benefit of an additional unit of a resource. The concept of marginal utility emerged during the marginal revolution in nineteenth-century economics \citep{marshall1890principles}.
%like L\'eon Walras, William Stanley Jevons, Carl Menger, Vilfredo Pareto and later Jules Dupuit and Alfred Marshall.
While costs are typically straightforward to assess, determining benefits can be more nuanced.
\textit{Willingness to pay} provides one avenue to estimate public benefits, which are often intangible and distributed across many stakeholders.
Willingness to pay centers on the amount someone would pay to have (or alternatively, accept) a particular intervention, and can be used to approximate the benefits of an intervention \citep{hanemann1991willingness}.
More generally, when the ratio of benefits to costs exceeds unity, benefits justify the financial outlay.
Underlying this rationale is Kaldor–Hicks efficiency, which states that the benefits to the public are sufficient to (hypothetically) compensate for the costs to society \citep{boardman2017cost}.

%DEPRECATE FOCUS ON SROI (it's in the nonprofit / fourth sector)
%Social return on investment (SROI) is the ``social impact of a business or nonprofit’s operations in dollar terms'' \citep{sroi}. Computing SROI is extremely context dependent. In the case of homeless services and shelters, standardized and accepted SROIs are nonexistent in the literature.  \cite{sroi}, provides 10 general guidelines to follow when determining SROI. Key amongst them are to include benefits as well as disbenefits; and to avoid double counting of particular benefits throughout the calculation.
%[MAYBE?? MV1 "Therefore, the marginal value of each dollar of donations spent in a country decreases in terms of the resulting beneficiaries’ welfare." [Keshvari Fard et al., MSOM 2022] ]

%PARAGRAPH 5 Decision criteria for comparing alternatives: incremental analysis
\looseness-1
A CBA typically compares a modest number of interventions, or alternatives, against the incumbent status quo to arrive at the most attractive alternative.
This analysis is carried out incrementally: the alternatives are placed in non-decreasing order of cost, and the incremental benefits and costs between the incumbent and the next-most expensive alternative are weighed using a BCR.
If the benefits of the increment outweigh the costs, the incumbent is replaced with the more attractive alternative; otherwise, the incumbent remains.
The analysis proceeds until all alternatives are exhausted.
In contrast to evaluating a small, finite set of user-defined alternatives, we propose a new algorithmic framework that implicitly represents evaluating combinatorially many alternatives to identify one that maximizes the BCR objective.
%Similar to incremental cost-effectiveness ratio used in health economics.
%First order all alternatives from lowest to the highest present worth of net costs (government costs less salvage value, if any)
%Sometimes values used for benefits, in the health economics literature, are costs-per-QALY of at least £20,000 to £30,000, could be up to £50,000 or even £100,000

\subsection{\emph{Fractional Programming}}\label{fracProg}
Fractional programming is an area of mathematical optimization concerned with optimizing a fractional, or ratio, objective function over a constrained system  \citep{fractionalprog1}. These optimization problems are of the form:
\begin{equation}\label{eq:fractionalprogram}
    \max \frac{N(x)}{D(x)} ,
\end{equation}
where $x \in \mathcal{D} \subseteq \mathbb{R}^n; N(x), D(x): \mathbb{R}^n \mapsto \mathbb{R}$; and $D(x) > 0 \ \forall \ x \in \mathcal{D}$. When $N(x), D(x)$ are affine and $\mathcal{D}$ is polyhedral, expression \eqref{eq:fractionalprogram} takes the form of a \textit{linear} fractional program:
\begin{equation}\label{eq:linearfracprog}
    \max\limits_{x \in \mathcal{D}} \left \{ \frac{c^{\top}x + \alpha}{d^{\top}x + \beta} \, \vert \, Ax \leq b \right \},
\end{equation}
\noindent where $c,d,x \in \mathbb{R}^n$;
$b \in \mathbb{R}^m$;
$\alpha, \beta \in \mathbb{R}$;
and $ A \in \mathbb{R}^{m \times n}$. A common assumption is positivity of the denominator, that is $d^{\top}x + \beta > 0 $. Fractional programming has a rich and diverse history since at least \cite{neumann1937},who first considered optimizing a fractional objective to find economic equilibria, considering a ratio between economic activity and the costs of such activity.
\cite{stancu2019ninth} provides a recent survey of a variety of fractional programming applications.

\looseness-1
While the solution to fractional programs of certain forms can be facilitated by transformation \citep{charnes}, the algorithm of \cite{Dinkelbach} can optimize a fractional objective function over a connected and compact domain while making far fewer assumptions on the structure of the objective function.
% equivalence to Newton's method here
\cite{Ibaraki} showed the equivalence of Dinkelbach's algorithm to Newton's method for this class of optimization problem, and as \cite{fractionalprogAlgos} referred to Dinkelbach's algorithm as a \textit{Newton method}, we will also follow this convention.
Newton's method \citep{Dinkelbach} can solve a fractional program of the form of \eqref{eq:fractionalprogram} by iteratively solving a sequence of related, linearized subproblems over the same domain.
\cite{schaibleSuperLinear} showed that Newton's method has superlinear convergence, which takes on increasing importance when a subproblem solve is computationally complex.
Recent technological and algorithmic developments have driven many advances in computational optimization, empowering increasingly larger optimization problems to be solved to global optimality \citep{bertsimas2016best} facilitating the solution of large-scale fractional programs.

While the BCR is paramount to public sector decision making, to date there are very few studies that advocate for the optimization of BCRs. Those studies that do exist cover secondary side applications in either the fractional programming methodology literature \citep{bradley1974fractional} or single applications that could only marginally be viewed as public sector \citep[see, e.g.][]{zappone2015energy, murdoch2007maximizing}. In \cite{liu2024} and \cite{Zhou2025}, the authors consider optimizing over a fractional objective to improve relief organization effectiveness in the disaster preparedness and humanitarian logistics context.
One recent exception is the work of \cite{parkGemma}, which considers optimizing a BCR to locate warehouses for food aid distribution by the World Food Program in Angola. Similar to our approach, the authors account for demand heterogeneity in a not-for-profit, supply constrained context.
Our study differs by introducing a general framework for optimizing BCR decisions that allows for a status quo scenario representing the current operational state and maximizes the BCR across many alternatives to assess marginal activity beyond the status quo.
Our case study weighs marginal resource allocation decisions concerning existing public sector operations, namely marginal allocations of additional NYC shelter and service capacity.
We are among the first to investigate the maximization of a BCR over a mixed-integer polyhedral set for public sector decision making, and the first to propose a general framework, which we now present.

\section{A New Algorithmic Framework for BCR Optimization}\label{Methods}
\looseness-1
Our algorithmic framework optimizes the benefit to cost ratio (BCR) and  %, which is subsequently demonstrated through a real-world case study in Section \ref{Model}.
rests upon the foundation of a complex, real-world decision context such as those found in the public sector.
We assume the mathematical representation of the problem encodes an existing state, the \textit{status quo}, together with decision variables that allow for movement away from the status quo.
While many decision-optimization contexts may yield prescriptions that represent significant departures from the status quo, practical considerations such as deployment and change management can complicate such prescriptions.
Due to high costs associated with substantial short-term operational modifications, such decision making may benefit from relative operational stability, both to mitigate risks and improve future flexibility.

\looseness-1
Traditional CBA assumes \textit{preidentified} actions through some generative process.
Our framework, in constrast, assumes no preidentifed actions, only constraints and (possibly) variables that express tangible \textit{action}, or \textit{movement}, beyond the status quo.
\textit{Action} constraints ensure some condition such as a minimum expectation or logical condition be satisfied, inducing alternatives implicitly as the constraints represent the feasible set of possible alternatives.
Action constraints may induce a vast array of alternatives, which our framework evaluates according to their associated benefits and costs, on a ratio-maximizing basis.
Our framework frees decision-makers from the labor of preidentifying alternatives, while still considering a wide range of feasible options, thereby identifying solutions 
beyond the status quo that are most cost-effective and socially desirable.
As decision makers dictate the required level of action, our framework is uniquely positioned to identify solutions that maintain relative operational stability while optimally moving beyond the status quo.

\looseness-1
We assume a decision-optimization problem $\mathcal{P}$ 
representing the system being optimized, composed of problem data $\Delta$ and constrained decision variables $x$ in a nonempty, bounded feasible region $\mathcal{D}$ formed by existing constraints. %, and a possible objective function $\mathcal{O}$.
%Without loss of generality, we assume the decision variables $x \in \mathbb{R}^n$ can be partitioned into decision variables $x^C \in \mathbb{R}^{n_1}$ that \textit{can} readily change, and decision variables $x^U \in \mathbb{R}^{n_2}$ that are \textit{unable} to readily change, with $n_1 + n_2 = n$.
Any feasible state of the system can be represented by an assignment of the variables $x$, and it is assumed there is opportunity for improvement through optimization.
Our algorithmic framework is presented in Procedure \ref{proc:framework}.

% ==============

% There may be existing efficiencies to be found in movement away from the status quo.
% The optimized status quo $x^{\star}$ is a state of the decision variables $x$ under $\Delta$ and $\mathcal{D}$.

% So when >= 1, marginal resources are permitted
% We can augment problem state with $\hat{S}$ and $\hat{x}$

% If $B(x,xhat) > C(x,xhat)$, the costs are not minded as they are in a ratio-optimizing way (somewhat indiscriminately) =>  as it does not indicate the sets of decisions which most effectively deploy available resources
% When $B(x,xhat)/C(x,xhat)$ is optimized, the costs are minded in ...
% As the actions dominate the profit (benefits less costs) optimization, any relationship between marginal cost and marginal benefit will be overly emphasized.
% In a rational relationship, there is more accountability between the marginal benefits and costs (not letting one term dominate the other).
% ==============

\begin{algorithm}

% \caption*{\textbf{Procedure 1 }The Algorithmic Framework}
\caption{}\label{proc:framework}
\caption*{\textbf{Procedure 1 }Algorithmic Framework for Benefit to Cost Ratio Optimization in PSDM}
\begin{enumerate}[leftmargin=*,itemsep=0pt,parsep=0pt, topsep=3pt, partopsep=3pt]
        \item Initialize decision problem $\mathcal{P}$: data $\Delta$, decision variables $x \in \mathbb{R}^n$, %status quo $\bar{x}$, 
        feasible region $\mathcal{D}$.
        \item Define marginal cost function $C : \mathbb{R}^n \to \mathbb{R}$ consisting of all relevant costs.
        \item Define marginal benefit function $B: \mathbb{R}^n \to \mathbb{R}$ capturing all benefits accrued to the public.
        \item Form $\frac{B}{C}$, the marginal BCR objective function to be maximized.
        \item Augment $\mathcal{D}$ with \textit{action} constraints $\hat{\mathcal{D}}$ (and possibly variables, $\hat{x}$) that ensure sufficient movement away from status quo under the BCR objective: $\mathcal{D} \leftarrow \mathcal{D} \cap \hat{\mathcal{D}}$ ($x \leftarrow x \cup \hat{x})$. 
        \item Optimize the BCR over $\mathcal{D}$ to identify maximizer $x^{\star}$ via appropriate means. 
    \end{enumerate}
\end{algorithm}

The algorithmic framework provided in Procedure \ref{proc:framework} identifies optimal marginal activity from the status quo.  %$\bar{x}$.
Data $\Delta$ encodes the structure of the decision making context including the status quo.  
The marginal costs and benefits of the system are expressed through $C: \mathbb{R}^n \to \mathbb{R}$ and $B: \mathbb{R}^n \to \mathbb{R}$, respectively.
The function $C$ captures marginal system costs, which can often be estimated through standard cost accounting procedures as they have concrete monetary values \citep{Gass-1994}. 
%Potential costs include expansion costs, new buildings, hiring and training of new staff, or campaign expenditures (e.g. an advertising campaign).
The benefits function $B$ measures the total marginal effect of decisions, or actions, upon public sector system stakeholders, and is typically more challenging to quantify.
Our \textit{action variables} represent the possible decisions to take, such as how many units of a service to provide.

\looseness-1
Optimizing the ratio $B/C$ first necessitates movement \textit{away from the status quo}, as the denominator 
gravitates toward attaining zero marginal cost, that is, the \textit{do nothing} option.
This movement is assured through a set of one or more \textit{action} constraints $\hat{\mathcal{D}}$ and a (possibly empty) set of one or more \textit{action} variables $\hat{x}$ that augment the original system constraints $\mathcal{D}$ and variables $x$ to ensure sufficient movement away from the status quo.
Optimization over a large number of decisions represented by action constraints and variables can reveal optimal marginal activity. 
In contrast, more conventional techniques, such as optimizing the difference between benefits less costs, cannot readily reveal marginal decisions, leading to a lack of insight in decision making when moving away from the status quo.
\cite{cseref2009incremental} also constrain movement from the status quo in the related study of \textit{incremental} optimization, while optimizing linear objectives of interest. Multiple algorithms that  optimize  a fractional objective function over a polyhedral set exist, each with unique properties, including rate of convergence and permitted functional forms \citep{fractionalprogAlgos}. For the case study investigated in Section \ref{Model}, we use Newton's method, outlined in Section \ref{Dinkelbach}.% Considering accuracy of the solution probabilistically might also be important to the problem context.

% This is done through considering multiple objective functions when optimizing over the polyhedral set $S$, as it allows the decision makers to work symbiotically with researchers to determine the best solutions moving forward.

%In Section \ref{Inferior Solutions}, we outline how this framework mathematically determines what solutions are inferior to others, while still gaining important insights from having seen them in the first place. 

% After determining the benefits function $B$ and cost function $C$, it is important to have a set bounded $S$ of feasible points for the benefits and costs being considered, which is likely to be defined through constraints to the variables in $B$ and $C$.

%Without a set of feasible points, optimization over the set $S$ is not well defined, as ratio optimization will have $C \to 0$, which implies that $\frac{B}{C} \to \infty$. Newton's method (\citep{Dinkelbach}) for optimizing a fractional objective also makes the added assumption that $C > 0, \forall x \in S$. Therefore, we need to make sure that $C > 0$ in order to utilize this powerful algorithm for the framework's optimization step. This could lead to there being no possible optimal solutions. Additionally, an unbounded set will make gaining insight from solutions generated difficult, as solutions may suggest decisions that are impossible to deploy in reality.

\subsection{\emph{Ratios as Generalized Differences}}\label{sec:ratiosasDiff}
% There are relationships between common objective functions such as benefit \textit{less} cost (profit) maximization and benefit \textit{to} cost ratio optimization as outlined in Procedure 1.
%We selected four objectives to demonstrate to decision-makers alternative view points through consultation with our partners in NYC. In the context of RHY, our partners deem the following objectives as important: $\max_{x \in S} B(x), \ \min_{x \in S} C(x), \ \max_{x \in S} B(x) - C(x), $ and $\max_{x \in S} \frac{B(x)}{C(x)}$. While these objective functions were developed with the particular context of RHY in mind, they have wide reach for a wide array of public section decision making problems. 

%\subsubsection{Ratios as Generalized Differences}
%From Newton's method (\citep{Dinkelbach}) (see Algorithm \ref{alg:Dinkelbach}), we see that an optimal value $q^*$ can optimize a fractional objective function $\frac{B}{C}$.
Consider generalized profit as expressed by the difference between benefits $B$ and costs $C$ for weights $\alpha, \beta \in \mathbb{R}_+$: 
\begin{equation}\label{eq:generalizedDifference}
    \max\limits_{x \in \mathcal{D}} \  \alpha B(x) - \beta C(x).
\end{equation}

\looseness-1
When $\beta = 0$, expression~\eqref{eq:generalizedDifference} reduces to $\max_{x \in \mathcal{D}} \  \alpha B(x)  = \alpha \max_{x \in \mathcal{D}} \   B(x)$, which shares the same set of optimal solutions as $\max_{x \in \mathcal{D}}  B(x) $. While likely impractical, it is a point of reference for maximized benefit without regarding  cost. In the absence of constraints to the contrary, maximizing this objective may incur significant cost and, in general, feature inefficient resource utilization.

Alternatively, setting $\alpha = 0$ in~\eqref{eq:generalizedDifference} yields $\max_{x \in \mathcal{D}} \  -\beta C(x)$, which shares the same set of optimal solutions as $\min_{x \in \mathcal{D}} \ C(x)$. Not regarding benefit, this approach represents a cost-minimal approach, taking no more action than necessary, while producing more efficient resource utilization.
%This allows decision makers to consider their budgetary constraints.
This style of objective function may be seen when revenues are fixed within a system, such as an organization completing work for a fixed level of income through a grant, contract, or appropriation.

Generalizing further, allowing $\alpha = \beta$ reduces expression~\eqref{eq:generalizedDifference} to $\max_{x \in \mathcal{D}} \  B(x) - C(x)$, or standard profit maximization.
As benefits and costs are equally weighted in this case, every action constraint and action variable now has an equal weight when considering the benefit that it brings to the optimization model. Because of this, optimizing this objective function will not reveal marginal decisions that provide the greatest improvement from the status quo.
Additionally, it is common in the public sector for costs to be expressed in monetary units while benefits are more challenging to express as they exist in widely differing units, thus making $\alpha = \beta$ highly unlikely.
%It is worth noting here that the units of measure for $B$ and $C$  can differ (e.g. Disability-Adjusted Life Years for $B$, and dollars for $C$). Furthermore, benefits are often nuanced and can be difficult to compute, such as health benefits, and benefits to a community. Additionally, it weighs the benefits equally with costs when computing the difference, making the assumption that $1$ benefit unit is equal to one $1$ cost unit.

\looseness-1
For any other values of $\alpha, \beta$, dividing through expression~\eqref{eq:generalizedDifference} by $\alpha$ gives a new coefficient $q = \beta / \alpha$ on the cost $C(x)$ without changing the resulting set of optimal solutions from maximizing $B(x)- qC(x)$ over $x \in \mathcal{D}$.
Strict profit maximization, that is, maximizing the benefit to cost \textit{difference}, can be seen when $q=1$.
Alternatively, maximizing the benefit to cost ratio $\frac{B(x)}{C(x)}$ over $x \in \mathcal{D}$ is tantamount to finding the $q^{\star} \in \mathbb{R}$ for which $\max_{x \in \mathcal{D}} \  B(x) - q^{\star}C(x) = 0$, as then $q^{\star} = \frac{B(x)}{C(x)}$, and any optimizer for one also optimizes the other \citep{Dinkelbach}.

While the choice of objective function is ultimately that of the decision maker, advantages to ratio optimization include that
%(1) the profit maximization form that arises when $q = 1$, that is, $\max_{x \in S} B - C$, would be unattractive to pursue if the corresponding optimizer $x^{\star}$ has a breakeven point when $B - C = 0$.
(i) it bypasses the somewhat arbitrary assumption equating $\alpha = 1$ unit of benefit to $\beta = 1$ unit of cost; (ii) the resulting BCR, when $B$ and $C$ share the same unit, is unitless; 
and
(iii) it tends to prescribe highly beneficial actions that are otherwise too removed from the cost-minimizing form $\min_{x \in \mathcal{D}}  C(x) $. Cost minimal optimization will ignore potential investments that may provide significant benefit for only a small cost.

\subsection{\emph{Contrasting Net-Benefit and BCR Optimization}} \label{whyProfitOptimizationFails}
The algorithmic framework ensures marginal activity away from the status quo by augmenting the system of original variables $x$ and constraints $\mathcal{D}$ with a set of action constraints $\hat{\mathcal{D}}$ and a (possibly empty) set of action variables $\hat{x}$.
This marginal activity incurs associated benefits and costs, and as such, it is instructive to consider these components in isolation, noting that not all need to assume a positive value.
The marginal costs may take the form of
%\begin{equation*}
$\hat{C}(x,\hat{x}) = \hat{C}_1(x) + \hat{C}_2(\hat{x}) + \hat{C}_3(x,\hat{x})$,
%\end{equation*}
and the marginal benefits as
%\begin{equation*}
$\hat{B}(x,\hat{x}) = \hat{B}_1(x) + \hat{B}_2(\hat{x}) + \hat{B}_3(x,\hat{x})$,
%\end{equation*}
yielding a marginal \textit{net benefit (benefits less costs)} function:
\begin{equation} \label{benefit_less_cost}
\hat{B}(x,\hat{x}) - \hat{C}(x,\hat{x}) = \hat{B}_1(x) + \hat{B}_2(\hat{x}) + \hat{B}_3(x,\hat{x}) - \hat{C}_1(x) - \hat{C}_2(\hat{x}) - \hat{C}_3(x,\hat{x}).
\end{equation}

For any $(x^{\star}, \hat{x}^{\star})$ that features $\hat{B}(x^{\star}, \hat{x}^{\star}) > \hat{C}(x^{\star}, \hat{x}^{\star})$ and is a maximizer for \eqref{benefit_less_cost}, optimization over $\mathcal{D}$ suggests activating combinations of variables that are collectively profitable, irrespective of the degree of profitability or corresponding cost for the individual variable.
In extremely resource-constrained settings, this may be ineffective as it is silent on the sets of decisions that most effectively deploy available resources.

Alternatively, the marginal BCR takes the form:
\begin{equation}
\frac{\hat{B}(x,\hat{x})}{\hat{C}(x,\hat{x})} = \frac{\hat{B}_1(x) + \hat{B}_2(\hat{x}) + \hat{B}_3(x,\hat{x})}{\hat{C}_1(x) + \hat{C}_2(\hat{x}) + \hat{C}_3(x,\hat{x})}.
\end{equation}

\looseness-1
Any optimizer $(x^{\star}, \hat{x^{\star}})$ of the marginal BCR is unlikely to exhibit the same behavior as an optimizer of \eqref{benefit_less_cost} due to the nonlinear relationship between the numerator and denominator; instead, there is more selectivity among activated variable combinations that improves the rational relationship between the marginal benefits and costs.
For example, if $\hat{B}(x, \hat{x}) = 10x + 5\hat{x}, \hat{C}(x, \hat{x}) = 5x + 3\hat{x}$ with binary $(x, \hat{x})$, optimizing $\hat{B}(x,\hat{x}) - \hat{C}(x,\hat{x})$ obtains an optimal solution of $(1,1)$ with cost $8$ and BCR $\frac{15}{8}$. By comparison, optimizing $\frac{\hat{B}(x,\hat{x})}{\hat{C}(x,\hat{x})}$ obtains an optima of $(1,0)$, reducing the cost to $5$ while improving the BCR to $2$. This is further demonstrated experimentally in Section \ref{CompuationalExperiments}.  

% \subsection{\emph{The Necessity of Action Constraints}}

% When optimizing benefit to cost ratios, it is necessary to constrain against the point $x$ which corresponds to ``doing nothing'' and not moving from the status quo. For a benefits function $f$ and a cost function $g$, the ratio must be constrained against the scenarios of minimal action. This allows the benefit to cost ratio to reveal the optimal \textit{marginal} decisions.
% \begin{theorem}
%     Let $f, g: \mathcal{D} \subseteq \mathbb{R}^n \to \mathbb{R}_{+}$ such that $g(x) > 0$ for all $x \in \mathcal{D}$. Further, assume that $f,g$ are both linear functions that are \textit{monotonically} increasing in each argument. This means that for $x \leq y$ (element wise), we have that $f(x) \leq f(y)$ and $g(x) \leq g(y)$. Further assume that 
%     \begin{subequations}
%         \begin{equation}
%             f(x) = a_f^{\top} x, g(x) = a^{\top}_g x , h(x) := \frac{f(x)}{g(x)}
%         \end{equation}
%         and that for all $x, y$ such that $x \leq y$, we have that 
%         \begin{equation}\label{eq:condition_thm}
%             \frac{a^{\top}_f x}{a^{\top}_g x} \geq \frac{a^{\top}_f (y-x)}{a^{\top}_g (y-x)}
%         \end{equation}
%     \end{subequations}
%     Then, we have that $h(x)$ is monotonically \textit{decreasing} in each argument. 
% \end{theorem}

% \begin{remark}
%     The condition \eqref{eq:condition_thm} can be interpreted as follows. The ratio of ``doing less'' is greater than the ratio of the marginal actions that took place between these two scenarios. 
% \end{remark}

\section{BCR Case Study: RHY Shelter and Service Capacity Expansion}\label{Model}
\looseness-1
To demonstrate our algorithmic framework, we focus on the problem of deploying marginal RHY shelter and serviced capacity in NYC, which has the largest homeless population in the US~\citep{usafacts_homelessness}. Additionally, the authors have a long-standing relationship with organizations and decision makers serving and working with RHY in NYC, providing a compelling case study for our analysis. Like many U.S. cities, New York is taking concrete steps toward a coordinated, system-level response to preventing and ending youth homelessness~\citep{Love2023}. We focus on transitional and independent living (TIL) shelters and programs as they offer housing and support services while youth work toward establishing independence. TIL programs empower youth by fostering growth in critical areas such as education, housing, employment, recreation, health, and safety, all of which promote self-sufficiency and independence \citep{Naccarato-2008}. 

%\subsection{\emph{Operations Research and Analytics in Human Services}}
%Operations research and analytics have seen rigorous use over the years to support decision-making in the public sector as well as nonprofit organizations, bringing substantial value and recognition \citep{Sinuany-2014}. Application areas include economic and operational problems related to public health, humanitarian and disaster aid, and public education \citep{besiou-2020, Rais-2011, Gass-1994}. 

\looseness-1

\subsection{\emph{Background: Reducing Human Trafficking Risk to RHY}}\label{s:lit_vulnerability}

RHY are at a particularly high risk of HT due to lack of access to basic necessities such as shelter, nutritious and affordable food, appropriate clothing, and healthcare \citep{Wright-2021}. Individuals and families encountering housing instability often experience long waits for public housing in the U.S. \citep{Gibbs-2014, arnosti}. Traffickers capitalize on these vulnerabilities to both entice RHY into trafficking and make it difficult to leave the exploitative situation  \citep{Middleton, martin2024identification}. One study on child sexual exploitation found that most victims in the study sample experienced homelessness or persistent housing instability \citep{dank}. Another study found that sex traffickers often spend time in areas known to be regular gathering spots for RHY and try to recruit young people by misrepresenting shelter availability and offering unsafe places to stay; similar tactics are known to be used by labor traffickers \citep{covenHouse}. Access to safe and stable housing is thus an effective way to reduce vulnerability to HT for RHY. %When housing and associated programmatic service needs remain under- or unaddressed, vulnerable individuals experience diminished welfare and are at increased risk of exploitation~\citep{TIP22}, calling for efficient deployment of public sector investments. 

\cite{Johnson-2012} advocate for a community-based operations research perspective that includes the interests of underrepresented, underserved, and vulnerable populations when designing solutions to improve economic efficiency, social equity, and administrative burdens.
Limited research has been conducted at the intersection of community-based operations research and scarce short-to-long-term housing resources to address homelessness \citep{Azizi-2018, Rahmattalabi-2022, Chan-2018, OriginalModel}.
Of course, housing instability extends beyond the homeless -- foster care systems face growing caseloads and high turnover in staff and foster families \citep{radey2022extent}, and refugee and asylum systems strain to meet rising needs of immigrants, frequently partnering with pro bono and community support to augment insufficient governmental support \citep{congress2022report}.

While operations research can be used to allocate scarce resources to address the increased HT risk of these youth due to housing capacity limitations, only a limited number of studies exist (\cite{MAASS2020100730,OriginalModel} \cite{kaya2022discrete}).

 \cite{MAASS2020100730} optimize the location of residential shelters serving trafficked persons in a manner that maximizes a measure of societal impact, while satisfying budget constraints. \cite{kaya2022discrete} use a discrete event simulation model to evaluate the necessary capacity expansion for a shelter that serves RHY, with particular attention to the needs of LGBTQ+ youth.

\cite{OriginalModel} propose a model that projects the cost-minimizing capacity to deploy in existing NYC RHY shelters to provide sufficient housing and support services. The authors consider stochastic youth arrivals and stay durations, services provided in a periodic fashion, and service delivery time windows. Two types of decision variables are used to model how a shelter can expand capacity: (i) extra in-house resource decision variables capture the amount of capacity that can be added within the existing facility; and (ii) overflow shelter decision variables that project the number of youth needing to be directed to a placeholder \textit{overflow} shelter when the within-shelter capacity cannot be further expanded. While there is utility in understanding the cost-minimizing capacity expansion to meet RHY needs in NYC through combining extra in-house resources and an overflow shelter, \cite{OriginalModel} leaves open the question of where, and in what quantities and priorities, these extra housing and support services should be allocated to appropriately match the needs of RHY. The authors focus on the cost of expanding services and, just as importantly, on the assignment of youth to shelters. Our case study to deploy marginal capacity to RHY shelters in NYC also features (i) extremely limited resources, (ii) multiple objectives, and (iii) uncertainty.
We address this resource scarcity and uncertainty by estimating the benefits and costs of capacity expansion alternatives, leveraging primary and secondary data. %, explained further in Section \ref{sec:datasets}.  %We extend \cite{OriginalModel} by allocating new shelter capacity in a manner that maximizes a societal BCR. Our approach considers how to strategically allocate limited resources across multiple candidate locations.

\subsection{\emph{Mathematical Model to Deploy Marginal Capacity}} \label{sec:marginalDeploymentModel}
\looseness-1
We extend the model of \cite{OriginalModel} to consider the marginal deployment of new shelter and service capacity in NYC. Rather than a single overflow location, we allow for more granular capacity expansion of the shelter system. We refer to capacity expansion as either adding additional short-term resources to existing organizations (i.e., adding capacity as needed each time period) or opening new shelter locations entirely. In addition to capacity expansion, we also model the assignment of youth to these organizations over time to meet their needs. The set of time periods we consider is represented by $T$.

Let $S$ be the set of shelter and service organizations for the model to consider. $S$ contains four disjoint subsets: $S=S^{sq-b} \cup S^{sq-r} \cup S^{new} \cup S^{service\ gap}$. $S^{sq-b}$ are the existing organizations that offer shelter beds and other services. $S^{sq-r}$ are the existing organizations that \textit{do not} offer shelter beds. However, organizations that do offer beds can \textit{refer} youth to organizations in $S^{sq-r}$ for services as needed. For example, a mental health clinic would be placed within $S^{sq-r}$. Here ``sq" indicates that the existing organizations are the current status quo. $S^{new}$ is the set of new organizations that can be opened. Assume that an organization $s \in S^{new}$ can only be opened if a minimum of $\mathcal{K}_{s}$ youth are assigned to it. 

The baseline capacity at organization $s\in S$ for service $i\in I$ at time $t\in T$ is $c_{s,i}^t$. However, there is flexibility to expand service capacity for one time period at a cost of $\gamma_{s,i}$ per unit increase. The maximum capacity that can be accommodated is $\mu_{s,i}$.

Each organization in $S$ has attributes related to its geographic location (e.g., NYC Borough), the youth demographics it serves, and the services it offers (e.g., beds, mental and physical health support, etc.). The set of services and demographics considered are denoted by $I$ and $D$, respectively. Binary variables $\sigma_{s,i}$ and $\beta_{s,d}$ taking value 1 indicate organization $s\in S$ offers service $i\in I$ and serves demographic $d\in D$, respectively . %The binary vector $\beta_s$ has length $|D|$ and is used to indicate whether organization $s\in S$ serves demographic feature $d\in D$ or not. 
$S^{service\ gap}$ is a set containing a single auxiliary organization to ensure model feasibility; youth placed in $S^{service\ gap}$ indicate that no other organization $S^{sq-b} \cup S^{sq-r} \cup S^{new}$ can satisfy the youth's needs and demographic.

 % The locations of new potential shelters are informed by \cite{OriginalModel} and determined through the algorithmic framework that selects the best marginal deployment of new shelter capacity,  opening new locations and expanding existing shelters. 
Let $Y$ be the set of RHY requiring placement. Each youth $y\in Y$ arrives at a discrete time period $l_y$ within the time horizon $T$. %, has a demographic profile $\alpha_y$ and requires a set of services from $I$. %The binary vector $\alpha_y$ has length $|D|$ and is used to indicate whether youth $y\in Y$ has demographic feature $d\in D$ or not. 
Youth $y\in Y$ requiring service $i\in I$ and being part of demographic $d\in D$ is indicated by binary variables $\nu_{y,i}$ and $\alpha_{y,d}$ taking the value 1, respectively. Youth must begin receiving their services between time periods $a_{y,i}$ (earliest start time) and $b_{y,i}$ (latest start time) and receive services for a duration of $d_{y,i}$ time periods. The parameter $f_{y,i}$ indicates the frequency at which youth must receive the service (e.g., continuous bed access, weekly mental health services).

Tables \ref{tab:modelSets}-\ref{tab:parameters1} list the sets and parameters in Appdendix \ref{mod:tables}. More detailed descriptions can be found in \citet[Table~2]{OriginalModel}.

\subsubsection{\emph{Decision Variables for Marginal Deployment}} \label{ss:dec_variables} 

The decision variables are summarized in Table \ref{tab:decisionVars}. Note that $\nu_s$ represents the action variables for this model.
% model has at most $|S^{new}| + |Y||S| + |Y||S||I| + |Y||S||I||T| + |S^{sq}||I||T|$ variables. 
For notational convenience, let $\boldsymbol{\nu}, \boldsymbol{\pi},\boldsymbol{U}, \boldsymbol{X}, \boldsymbol{E}$ denote the vectorization of each set of decision variables. 
\begin{table}[H]
\centering 

\renewcommand{\arraystretch}{0.7} % Decrease to compress rows
\begin{tabular}{lcl}
\toprule 
Variable        & Domain             & Definition   \\ \midrule
$\nu_s$         & $\{0,1 \}$           & $1$ if (new) organization $s \in S^{new}$ is opened\\
$\pi_{y,s}$ & $\{0,1\}$ & $1$ if youth $y \in Y$ is ever assigned to organization $s \in S$\\
$U_{y,s,i}$     &  $\{0,1\}$      & $1$  if youth $y \in Y$ is ever assigned to organization $s \in S$ for service $i \in I$  \\ 
$X^{t}_{y,s,i}$ & $\{ 0,1 \} $           & $1$ if youth $y \in Y$ is assigned to organization $s \in S$ for $i \in I$ at time $t \in T$     \\ 
$E^{t}_{s,i}$   & $\ZZ_{+} $ & Integer units of service $i \in I$ to add to organization $s \in S^{sq}$ at time $t \in T$ \\
\bottomrule
\end{tabular}
\caption{Decision variables used in the marginal deployment optimization model.}
\label{tab:decisionVars}
\end{table}

\subsubsection{\emph{Objective Function for Marginal Deployment}} \label{ss:objective}
The objective function consists of a benefit and cost function. We first discuss the benefit function.

\noindent \textbf{\underline{Benefit Function}} We model the benefits function as:
\begin{equation}\label{eq:benefits}
    B(\boldsymbol{U}, \boldsymbol{\nu}) = \sum_{y \in Y} \sum_{s \in S}  \sum_{i \in I} F_{y,s,i}  f_{y,i} U_{y,s,i}  + \sum_{s \in S^{new}}  \tilde{p}_{s} \nu_s,
\end{equation}
with the first term representing the benefits per youth served and the second representing benefits from each new organization opened. 

\looseness-1

%\subsubsection{Youth Benefits}\label{sec:YouthBenefits}

% The cost figure $p_{s}$ approximates system benefits through \textit{estimation by analogy}, a common technique for complex estimation circumstances that can provide reasonable values~\citep{karlsen2005management,shepperd1997estimating}.
% To estimate a reasonable public sector return on investment, we use estimates obtained from \cite{sroiValue}, who evaluated the social return on investment from opening a homeless shelter in London, UK through interviews with past and present service users of the shelter.
% With this multiplier $\rho$, the returns $\tilde{p}_{s}$ from opening an organization are given as $\tilde{p}_{s} = \rho \times p_{s}$.
% \begin{equation}\label{eq:sroiValue}
%     \tilde{p}_{s, \ell} = \rho \times p_{s,\ell}.
% \end{equation}

\noindent \textbf{\underline{Cost Function}}
Our total cost function $C$ is given by: 
\begin{equation}\label{eq:costs}
    C(\boldsymbol{\nu},\boldsymbol{E}, \boldsymbol{U}) = \sum_{s \in S^{new}} p_{s} \nu_s +   \sum_{s \in S^{sq} } \sum_{i \in I}\sum_{t \in T}  \gamma_{s,i} E^t_{s,i} + \sum_{y \in Y} \sum_{s \in S} \sum_{i \in I} r_{y,s,i} f_{y,i} U_{y,s,i} 
\end{equation}
The first component of the total cost function is the cost of opening a new organization $\nu_s$. The second  concerns \textit{expansion} of available services within organizations. For an organization $s \in S^{sq}$, service $i \in I$ can be expanded within limits for a cost $\gamma_{s,i}$ per unit. The final component sums the cost of assigning youth $y \in Y$ to organization $s \in S$ for service $i \in I$. This component has assignment costs $r_{y,s,i}$ that reflects the cost to provide services to youth differs by organization.

With these components, we formulate the BCR objective function as:
\begin{equation}\label{eq:ratioObjectiveFunciton}
    \frac{B(\boldsymbol{U}, \boldsymbol{\nu})}{C(\boldsymbol{\nu},\boldsymbol{E}, \boldsymbol{U})} = \frac{\sum\limits_{y \in Y} \sum\limits_{s \in S}  \sum\limits_{i \in I} F_{y,s,i}  f_{y,i} U_{y,s,i}  + \sum\limits_{s \in S^{new}} \tilde{p}_{s} \nu_s}{\sum\limits_{s \in S^{new}} p_{s} \nu_s +   \sum\limits_{s \in S^{sq} } \sum\limits_{i \in I}\sum\limits_{t \in T}  \gamma_{s,i} E^t_{s,i} + \sum\limits_{y \in Y} \sum\limits_{s \in S} \sum\limits_{i \in I} r_{y,s,i} f_{y,i} U_{y,s,i} }
\end{equation}

\subsubsection{\emph{Action Constraints for Marginal Deployment}} \label{ss:constraints}
In light of the fractional objective function that naturally drives the denominator toward zero, we augment the constraints proposed in \cite{OriginalModel} with new constraint sets to promote action beyond the status quo:

\vspace*{-10mm} % Adjust the negative value as needed
\looseness-1
\begin{subequations}\label{eq:constraints}
\begin{align}
& & & \sum\limits_{s \in S^{new}} \nu_{s} \geq \lambda, &  
\label{eq:AtLeastlambdaOpenShelters} \\
& & & \sum\limits_{y \in Y} \pi_{y,s} \geq \mathcal{K}_{s} \nu_s, & \forall s \in S^{new} \label{eq:CriticalMass} \\ 
& & & \sum\limits_{y \in Y} \sum\limits_{i \in I} U_{y,s,i} \leq |I| |Y| \nu_s, & \forall s \in S^{new} \label{eq:newShelterOpen1}\\
& & & \sum\limits_{i \in I} U_{y,s,i} \leq |I| \pi_{y,s}. & \forall y \in Y, \forall s \in S \label{eq:NewShelterOpen2} 
\end{align}
\end{subequations}

Constraint \eqref{eq:AtLeastlambdaOpenShelters} ensures we open at least $\lambda$ new shelter programs. 
Constraint set~\eqref{eq:CriticalMass} ensures that, for any new shelter program $s$ to open, a critical mass of at least $\mathcal{K}_{s}$ youth must be assigned to the program.
Constraint set~\eqref{eq:newShelterOpen1} ensures that youth may be assigned to services in a program only when it is open, while constraint set~\eqref{eq:NewShelterOpen2} allows for support services to be provided to youth who are assigned to specific organizations.
%Constraint set~\eqref{eq:constraints} admits at most $|Y||S| + 2|S^{new}| + 1$ constraints.
\looseness-1

% \newpage

\begin{subequations}
\label{eq:constraints}

% \noindent\textit{Added short-term capacity at existing organizations cannot exceed maximum limit:}

\vspace{-5mm}

\begin{equation}
\hspace*{15mm} c_{s,i}^t + E^{t}_{s,i} \leq \mu_{s,i},
\quad \forall s \in S^{sq},\ i \in I,\ t \in T.
\label{con:shelter_limit}
\end{equation}

% \vspace{2mm}
% \noindent\textit{Assigned youth cannot exceed capacity at existing organizations:}

\vspace{-5mm}

\begin{equation}
\hspace*{15mm} \sum_{y \in Y} X^{t}_{y,s,i} \leq c_{s,i}^t + E^{t}_{s,i},
\quad \forall s \in S^{sq},\ i \in I,\ t \in T.
\label{con:capacity_existing}
\end{equation}

% \vspace{2mm}
% \noindent\textit{Assigned youth cannot exceed capacity at new organizations:}

\vspace{-5mm}

\begin{equation}
\hspace*{15mm} \sum_{y \in Y} X^{t}_{y,s,i} \leq c_{s,i}^t,
\quad \forall s \in S^{new},\ i \in I,\ t \in T.
\label{con:capacity_new}
\end{equation}

% \vspace{2mm}
% \noindent\textit{Activate assignment indicator $U_{y,s,i}$ if youth is ever recieves the service from the organization:}

\vspace{-5mm}

\begin{equation}
\sum_{t \in T} X^{t}_{y,s,i} \leq |T|\,U_{y,s,i},
\quad \forall y \in Y,\ s \in S,\ i \in I.
\label{con:U_active}
\end{equation}

% \vspace{2mm}
% \noindent\textit{Avoid assignments outside of the permissible start service time window:}

\vspace{-5mm}

\begin{equation}
\sum_{s \in S}\sum_{t=0}^{a_{y,i}-1} X^t_{y,s,i}
\;+\;
\sum_{s \in S}\sum_{t=b_{y,i}+d_{y,i}+1}^{|T|} X^t_{y,s,i}
= 0,
\quad \forall y \in Y,\ i \in I.
\label{con:not_assigned_outside_range}
\end{equation}

% \vspace{2mm}
% \noindent\textit{Provide service assignment within the youth's permissible start time window:}

\vspace{-5mm}

\begin{equation}
\sum_{s \in S}\sum_{t=a_{y,i}}^{b_{y,i}} X^t_{y,s,i} \geq 1,
\quad \forall y \in Y,\ i \in I.
\label{con:assigned_at_least_once}
\end{equation}

% \vspace{2mm}
% \noindent\textit{Receive the service at the needed frequency during the needed duration:}

\vspace{-5mm}

\begin{equation}
\sum_{s \in S}\sum_{t=0}^{b_{y,i}+d_{y,i}} X^t_{y,s,i} = f_{y,i},
\quad \forall y \in Y,\ i \in I.
\label{con:assigned_sufficiently}
\end{equation}

% \vspace{2mm}
% \noindent\textit{Only assign youth to programs that provide their needed services:}

\vspace{-5mm}

\begin{equation}
%\sum_{i \in I}\sum_{t \in T} X^t_{y,s,i} = 0,
\sum_{t \in T} X^t_{y,s,i} = 0,
\quad \forall i\in I, y \in \{Y|\eta_{y,i}=1\},\ s \in \{S|\sigma_{s,i}=0\}
\label{con:service_match}
\end{equation}

% \vspace{2mm}
% \noindent\textit{Program must serve the youth's demographics:}

\vspace{-5mm}

\begin{equation}
\sum_{i \in I}\sum_{t \in T} X^t_{y,s,i} = 0,
%\quad \forall y \in Y,\ s \in S \text{ where } \exists\, n \text{ s.t. } (\beta_{s}-\alpha_{y})_{n}<0.
\quad \forall d\in D, y \in \{Y|\alpha_{y,d}=1\},\ s \in \{S|\beta_{s,d}=0\}
\label{con:demographic_match}
\end{equation}

% \vspace{2mm}
% \noindent\textit{A youth may receive each service from at most one organization in each time period:}

\vspace{-5mm}

\begin{equation}
\sum_{s \in S} X^{t}_{y,s,i} \leq 1,
\quad \forall y \in Y,\ i \in I,\ t \in T.
\label{con:in_one_place_at_a_time}
\end{equation}
\end{subequations}

% %
\noindent The remaining constraint sets in the mathematical model inherit from~\cite{OriginalModel}, with slight notational modifications\footnotemark \footnotetext{Constraint sets $(2c), (4b),$ and $(4c)$ are not inherited from the model, because we allow youth to go to multiple shelters for the same service, and we do not model periodic services.}. 
Constraint system \eqref{con:shelter_limit} limits the maximum amount that can be expanded in existing organizations. Constraint systems \eqref{con:capacity_existing} and \eqref{con:capacity_new} enforce capacity limits for organizations during youth assignment.  Constraint system \eqref{con:U_active} ensures that $U_{y,s,i} = 1$ if a youth is assigned to an organization for a service at any time $t \in T$. Constraint system \eqref{con:not_assigned_outside_range} prevents youth assignment beyond their available start time window. Constraint system \eqref{con:assigned_at_least_once} requires a youth to be assigned in their start time window, while set \eqref{con:assigned_sufficiently} guarantees they receive the requested service amount. Constraint systems \eqref{con:service_match} and \eqref{con:demographic_match} restrict youth assignment to service and demographic compatible organizations. Finally, constraint system \eqref{con:in_one_place_at_a_time} ensures each youth is assigned to at most one place at a given time for each service. In Appendix \ref{app:presolve} we include standard cuts that do not impact the meaning of the model but improve solvability.

\subsection{\emph{Newton's Method}} \label{Dinkelbach}
We employ Newton's method \citep{Dinkelbach} to maximize the BCR as expressed in \eqref{eq:ratioObjectiveFunciton}, with pseudocode provided in Algorithm \ref{alg:Dinkelbach}.
Newton's method has superlinear convergence, limiting the number of iterations prior to convergence \citep{schaibleSuperLinear}. 
Newton's method requires two continuous functions $B$ and $C$ and further assumes that $C$ is positive. Action constraints as outlined in Section \ref{ss:constraints} ensure that $C(x) > 0$ for all $x \in \mathcal{D}$. We assume that feasible region $\mathcal{D}$ is bounded and nonempty. The compactness of $\mathcal{D}$, the positivity of $C$, and the continuity of $B$ and $C$ (when applying a linear relaxation) suffice to guarantee existence of, and superlinear convergence to, a global optimizer.
\setcounter{algorithm}{0}
\begin{algorithm}[H]
\caption{Newton's Method for Maximizing the Benefit to Cost Ratio.} \label{alg:Dinkelbach}

\begin{enumerate}[leftmargin=*, itemsep=0pt, parsep=0pt, topsep=3pt, partopsep=3pt]
    \item Initialize $k=0$ and $q_k = 0$, and tolerance $\varepsilon >0$.
    \item Solve $F(q_k) = \max \{ B(x) - q_k C(x) \ | \ x \in \mathcal{D}\}$, and let $x_k^{\star}$ be its solution.
    \item If $F(q_k) < \varepsilon$: Stop; $x_k^{\star}$ is the optimizer and $q^{\star} = q_k^{\star}$ is the optimal ratio.
    \item Else: Update $q_{k+1} = \frac{B(x_k^{\star})}{C(x_k^{\star})}$, $k=k+1$, and go to Step 2.
\end{enumerate}
\end{algorithm}
% \begin{algorithm}[H]
% \caption{Newton's Method for Maximizing the Benefit to Cost Ratio.} \label{alg:Dinkelbach}
% 1. Initialize $k=0$ and $q_k = 0$. \newline 
% 2. Solve $\ \ F(q_k) = \max \{ B(x) - q_k C(x) \ | \ x \in \mathcal{D}, \ $ and let $x_k^{\star}$ be its solution. \newline 
% 3. If $F(q_k) < \varepsilon:$ Stop; $x_k^{\star}$ is the optimal solution and $q^{\star} = q_k^{\star}$ is the optimal ratio. \newline 
% 4. If $F(q_k) \geq \varepsilon: $ Update $q_{k+1} = \frac{B(x_k^{\star})}{C(x_k^{\star})}$, $k=k+1$, and go to Step 2.
% \end{algorithm}

\subsection{CORE: A Simple Heuristic for Large Instances}\label{sec:core}

We introduce a batch heuristic, \textit{Commit Openings, Refresh, and Expand} (CORE), to improve performance of our approach under larger youth counts. CORE operates by solving for a batch of randomly selected youth from $Y$, committing the non-youth-specific decisions from that solution (i.e., program openings $\nu_s$ and expansion allocations $E_{s,i}^t$), updating parameters (e.g. capacities $c_{s,i}^t$), and then re-solving with the next batch, until all youth are assigned. Algorithm~\ref{alg:core} outlines how CORE operates.

\begin{algorithm}[H]
    \caption{Commit Openings, Refresh, Expand (CORE)}\label{alg:core}

    \begin{algorithmic}
       
        \STATE Initialize $Y \leftarrow Y^0$, $\ell \leftarrow 1$
\REPEAT
    \STATE 1. Select a random subset of youth $Y^{\ell} \subset Y$ uniformly without replacement
    \STATE 2. Solve model on $Y^{\ell}$
    \STATE 3. Update $Y \leftarrow Y \setminus Y^{\ell}$; update capacities and related parameters (e.g. $c_{s,i}^t$); commit non-youth-specific decisions from this solution, namely program opening variables $\nu_s$ and expansion variables $E_{s,i}^t$
    \STATE $\ell \leftarrow \ell + 1$
\UNTIL{all youth in $Y^0$ are assigned}
    \end{algorithmic}
\end{algorithm}

% \begin{enumerate}[leftmargin=*, itemsep=0pt, parsep=0pt, topsep=3pt, partopsep=3pt]
    
%     \item Solve the model on a randomly chosen subset of youth $\bar{Y} \subset Y$; fix the non-youth-specific decisions, from this solution, namely the program opening variables $\nu_s$ and expansion variables $E_{s,i}^t$.
%     \item Update capacities and related parameters (e.g. $c_{s,i,t}$), add the next batch of youth, and re-solve, using already opened shelters $\nu_{s}$ before considering new openings.
    
% \end{enumerate}

\section{Computational Studies}\label{CompuationalExperiments}
\looseness-1
We evaluate the proposed techniques for marginal expansion through computational experiments using large, stakeholder-informed datasets that are representative of the real-world. 
The experiments are designed to assess three questions: how a BCR objective differs from conventional cost and net-benefit objectives; how required action levels affect marginal expansion decisions; and whether the CORE heuristic preserves high-quality recommendations at larger youth counts.

\noindent \textbf{Computational Setup.}
\looseness-1
All experiments were implemented in Python and solved with Gurobi Optimization 11.0.2 \citep{gurobi} on Northeastern University's Explorer Research Cluster. For all experiments, we set the Gurobi NodeFileStart parameter to 0.5 and used a 5\% relative MIPGap tolerance.
Limited testing with smaller gap values (e.g., 1\%) did not improve objective values beyond the reported solutions, while significantly increasing runtimes. Each instance received up to the maximum allowed 48 hours of wall clock compute time.

\looseness-1

% Datasets 
\subsection{\emph{Datasets}}\label{sec:datasets}
% We model weekly capacity allocation decisions over a $6$ month horizon, with time periods, service duration and frequency specified in weeks.
% In this section, we detail the datasets used in our computational experiments.
\subsubsection{\emph{RHY Organization and Youth Profiles}}

This study uses three datasets: i) demographic and capacity profiles for existing RHY shelter programs; ii) demographic and capacity profiles for potential new RHY shelter programs; and iii) demographic and service need profiles for RHY youth. These datasets define the base computational setting: $|S^{sq}| = 8$ existing shelter programs, $|S^{new}| = 30$ potential new shelter programs, $|Y| = 168$ youth in exclusively Transitional Independent Living (TIL) programs, $|I| = 7$ support services, and $|T| = 24$ weekly time periods. 
The 8 existing shelter program profiles are drawn from~\cite{OriginalModel}, and additional details on the set of potential new shelters considered are provided in Appendix~\ref{dataDetailsAppendix}.
The seven support services are beds; health, including mental health, physical health, and substance use support; personal development, including education, employment, and life and practical skills; housing search, including long-term and crisis; legal support; parental support; and financial support.  

For $s \in S^{sq-b}$, service availability $\sigma_{s,i}$ and demographic availability $\beta_{s}$ are determined through interviews with shelter partners.
For referral programs $s \in S^{sq-r}$, we assume demographic eligibility is unrestricted and that  each program focuses on a single service $i \in I$.
This construction yields $|S^{sq-r}| = 6$ referral programs.
For potential new shelter programs $s \in S^{new}$, we define six demographic eligibility profiles $\beta_{s,d}$: baseline\footnotemark \footnotetext{Baseline refers to a majority vote of what each shelter $s \in S^{sq-b}$ contains. Thus, it can be considered an ``average'' shelter.}, male only, female only, 16-21 only, 22-23 only, and one that explicitly serves RHY parenting youth. Further details are provided in Appendix~\ref{dataDetailsAppendix}). For each service $i \in I$, the capacity of a potential new shelter is set equal to the corresponding average capacity across the existing shelter programs $s \in S^{sq-b}$.  
\looseness-1

In the experiments, these base profiles remain fixed, while selected operational parameters are sampled to construct multiple instances.
Specifically, we vary youth arrival times, shelter capacities, and costs to examine sensitivity to plausible variation in demand timing, resource availability, and expansion costs.

\subsubsection{\emph{Model Data}}\label{sec:costs}
\looseness-1
Between December 2021 and June 2023, the authors and a group of social scientists [blinded for peer review] 
%\citep{dank2024understanding}
surveyed 384 RHY in NYC to collect information on housing and service needs. We use these data to generate youth request parameters $\eta_{y,i}$ and demographic profiles $\alpha_{y}$. To assess solution sensitivity, we generate $10$ sampled instances that vary youth arrivals, shelter capacities, and costs while holding the underlying shelter and youth profiles fixed. Appendix \ref{dataDetailsAppendix} provides additional details on the data generation procedure.

\textbf{Youth Arrival and Service Timing.} We model the youth arrival time as $l_y \sim U\{1,...,|T|\}$, where $U\{\cdot\}$ denotes the discrete uniform distribution, so that youth requests may occur uniformly throughout the planning horizon. To represent limited flexibility in when youth can receive a given service, we set the earliest start time for service $i$ to $a_{y,i} = l_y$ and the latest acceptable start time to $b_{y,i} = a_{y,i}+\Delta, \ \Delta \sim U\{0,1,2\}$. This construction represents the small but realistic scheduling window that case managers report when placing youth into services. Youth then remain in the system for $d_{y,i} \sim U\{1,...,|T|-1\}$ weeks and receive service $i$ a total of $f_{y,i}$ times; for beds, we set $d_{y,i}=f_{y,i}$ to represent a continuous stay. These values are based on provider-reported lengths of stay, which typically range from one week to nearly the full planning horizon of $6$ months.

\textbf{Shelter Capacities.} Each shelter $s \in S^{sq}$ has a time-indexed service capacity $c_{s,i,t}$ drawn from a triangular distribution informed by interviews with our NYC shelter program partners. Shelter staff also report informal yet common ways to extend capacity when shelters are full, such as counselors staying overtime or providing hotel vouchers. To represent this operational flexibility, we define a maximum capacity $\mu_{s,i}$ set to $1.3$ times the corresponding shelter service capacity, $c_{s,i,t}$.
For referral organizations $s \in S^{sq-r}$, we set $c_{s,i,t}$ to a large value for each relevant service $i \in I$ and time $t \in T$ representing effectively unconstrained referral capacity.
For the remaining existing shelters, we scale $c_{s,i,t}$ to represent only the remaining $5-10\%$ available capacity, consistent with partner reports that most shelters operate at $90-95\%$ occupancy and rarely have every bed available. We model the existing shelters being nearly full as this is the reality of what is on the ground. The youth that we model in the system are those that will inform our expansion decisions, while the existing shelters continue to serve others who also have needs. 
For each new shelter $s \in S^{new}$, we set each service capacity equal to the average capacity across existing organizations $s \in S^{sq}$. For each potential new shelter $s \in S^{new}$, we set $\mathcal{K}_{s}$=1, so opening the shelter requires at least one assigned youth.

\textbf{Cost Estimation.} We estimate $\gamma_{s,i}$ using public data, including hourly salary estimates and NYC hotel voucher prices \citep{Gross-2021, NYCHRA-2022, Ziprecruiter-2022}. We generate the cost $p_{s}$ of opening a new shelter using a 12-bed capacity assumption and a cost per bed of $C_{\text{bed}} =\$57,778$. We compute borough multipliers from~\cite{furmanData} average 2022 rent values, scaling each borough's values relative to the mean. 
Table \ref{tab:BoroughMult} reports the resulting shelter cost multipliers for the five boroughs of NYC.

\begin{table}[H]
\centering 

\renewcommand{\arraystretch}{0.8} % Decrease to compress rows
\begin{tabular}{lccccc}\toprule 
               % & The Bronx & Staten Island & Queens & Brooklyn & Manhattan \\ \midrule % 1 = Bronx, 2 = Staten Island, 3 = Queeens, 4 = Brooklyn, 5 = Manhattan
             Borough  & 1 & 2 & 3 & 4 & 5 \\ \midrule 
Cost Multiplier $L_{s}$& 0.792     & 0.980         & 1.048    & 0.991     & 1.189     \\ \bottomrule 
\end{tabular}
\caption{The shelter cost multipliers for the five boroughs of NYC.}
\label{tab:BoroughMult}
\end{table}

For shelter programs $s \in S^{sq} \cup S^{new}$, we normalize the assignment cost $r_{y,s,i}$ to $1$.
For referral programs $s \in S^{sq-r}$, we generate assignment costs $r_{y,s,i}$ from a triangular distribution calibrated to make in-house placement preferable to referral assignment when feasible. The distributional parameters are informed by input from our partners in NYC, with additional details provided in Appendix ~\ref{dataDetailsAppendix}. 
Referral programs $s \in S^{sq-r}$ are assigned substantially larger capacities than shelter programs $s \in S^{sq} \cup S^{new}$, yet they do not provide bed services.
%We change the value to $1$  to guarantee $C(X^t_{y,s,i}, E^t_{s,i}, \nu_s) > 0$ to satisfy the assumptions of Newton's method \citep{Dinkelbach}. 

\textbf{Benefits Estimation.} Quantifying $F_{y,s,i}$, the benefit from serving youth $y$ at organization $s$ for service $i$, is challenging.
We assign zero incremental objective benefit to placements in referral programs $s \in S^{sq-r}$.
This assumption is not a statement about the intrinsic value of referral services; rather, because they are modeled as having effectively unlimited capacity, the optimization focuses on decisions that expand constrained in-house services.
For existing shelter programs and prospective new organizations, a natural alternative would be to adopt a willingness-to-pay framework similar to \cite{MAASS2020100730}. However, because RHY typically lack the financial means to pay for services, willingness to pay would understate service value and may assign zero value to services that are substantively appropriate. 

Quantifying $\tilde{p}_{s}$, the benefit associated with opening a new shelter organization $s \in S^{new}$, is also challenging: new shelters impose costs on city governments, yet they may generate substantial benefits for youth, neighborhoods, and surrounding communities.

Because these benefit parameters are difficult to estimate directly, we conduct sensitivity analysis over alternative parameterizations.
We consider two parameterizations: one assigns raw values directly to each benefit coefficient, while the other scales costs by a multiplier to represent a return-on-investment interpretation.
Under the multiplier interpretation\footnotemark \, youth service benefits are scaled from assignment costs as $F_{y,s,i} = m \times r_{y,s,i}$  \footnotetext{Here, it is also beneficial for us to have the benefit for a referral being $0$. If it were greater than $0$, we would see solutions that actually \textit{incentivize referrals over in-house assignment} due to the fact that referrals carry a higher cost. Future work can investigate more granular multipliers than what is considered here.} and opening benefits are scaled from shelter-opening costs as $\tilde{p}_{s} = m' \times p_{s}$ ;under the raw value parameterization, benefit coefficients are assigned directly. Table \ref{tab:benefit-values} lists the raw benefit values and multipliers used in the sensitivity analysis. 

\begin{table}[h]
  \centering
  \renewcommand{\arraystretch}{1.4}
  \begin{tabular}{l | p{5.3cm} | p{5.3cm}}
  \toprule
   & Youth & Shelter \\
  \midrule
  Raw values &
    \makecell[l]{ $F_{y,s,i} \in \{ 1,\ 10,\ 100,\ 500,$\\$1000,\ 5000,\ 10000 \} $} &
    \makecell[l]{$\tilde{p}_{s} \in \{ 10,\ 100,\ 500,\ 1000,$\\$5000,\ 10000 \}$} \\
  \midrule
  Multiplier &
    \makecell[l]{$m \in \{  1.0,\ 1.5,\ 2.0,\ 2.5,\ 3.0,\ 3.5,$\\$4.0,\ 5,\ 15,\ 25,\ 35,\ 45 \}$} &
    \makecell[l]{$m' \in \{ 0.0001,\ 0.0005,\ 0.001,\ 0.005,$\\$0.01,\ 0.05,\ 1.0,\ 1.5,\ 2.0,$\\$2.5,\ 3.0,\ 3.5,\ 4.0,\ 5,\ 10,$\\$15,\ 20,\ 25,\ 30,\ 35,\ 40,\ 45,\ 50\}$} \\
  \bottomrule
  \end{tabular}
  \caption{Raw benefit values and multipliers for youth and shelter benefits considered in sensitivity analysis.}
  \label{tab:benefit-values}
\end{table}
\subsection{\emph{Computational Experiments}}\label{sec:ComputationalExperiments}

\looseness-1
We conduct three experiments to evaluate how the BCR objective guides marginal resource deployment.
These experiments evaluate BCR optimization through an objective-function comparison, a minimum shelter count sensitivity analysis, and an adaptive heuristic scalability test. % for optimal marginal resource deployment.

First, we study capacity expansion under the BCR objective  $\max\limits_{x \in \mathcal{D}}  \frac{B(x)}{C(x)}$ and, akin to Section \ref{sec:ratiosasDiff}, compare its recommendations with conventional cost minimization $ \min_{x \in \mathcal{D}}C(x)$ and net-benefit objectives $\max_{x \in \mathcal{D}}B(x) - C(x)$.
Second, we vary the minimum number of shelter programs required to open across NYC by setting $\lambda \in \{1,2,3,4,5,6\}$ in action constraint \eqref{eq:AtLeastlambdaOpenShelters}, allowing us to assess how required expansion levels affect the BCR.
Third, we evaluate scalability by varying the number of youth $|Y|$ and testing a heuristic designed to handle larger youth populations.

While ideally $X^{t}_{y,s,i}$ would be binary, we relax $X^{t}_{y,s,i}$ to be continuous and, in limited computational testing, find that over $99\%$ of these variables assume binary values.
The model emphasizes shelter-opening decisions $\boldsymbol{\nu}$ rather than the assignment-timing variables $\boldsymbol{X}$.
The model estimates which marginal deployment actions are needed to improve upon the status quo. 
\looseness-1
% Referrals are modeled as an absorption state for unmet demand, and are assumed to have infinite capacity. 
We  evaluate and compare solutions using \textit{system  utilization}, defined as: 
\begin{equation}
    \frac{\textit{Total number of youth service requests fulfilled in-house}}{\textit{Total youth number of service requests fulfilled}}. 
\end{equation}
\looseness-1
Because referral organizations are modeled with very high capacities, this metric captures the share of fulfilled service requests assigned to in-house shelter programs rather than referral organizations.
A high in-house fulfillment share indicates that many youth receive necessary services in-house, a desirable outcome that supports continuity of care, builds trust, and reduces transportation needs.

We also examine \textit{program occupancy} for new shelter programs, defined as the average fraction of program capacity used across all time periods. Specifically, program occupancy is defined as the average share of program capacity utilized over the planning horizon: 
\begin{equation}
    \frac{1}{|T|}\sum_{t = 1}^{T}\frac{\textit{Total service requests fulfilled by program at time }t}{\textit{Total capacity of program at time } t}.
\end{equation}

Together, these two performance measures capture both the extent to which youth are served in-house and the intensity with which newly opened shelter capacity is used.

\subsection{\emph{Model Results}}\label{baselineResults}
\looseness-1

We next present the model results, beginning with performance under alternative objective functions.

\subsubsection{\emph{Experiment 1: Performance under Different Objectives}}

\begin{figure}[H]
    \centering
    \includegraphics[width=\linewidth]{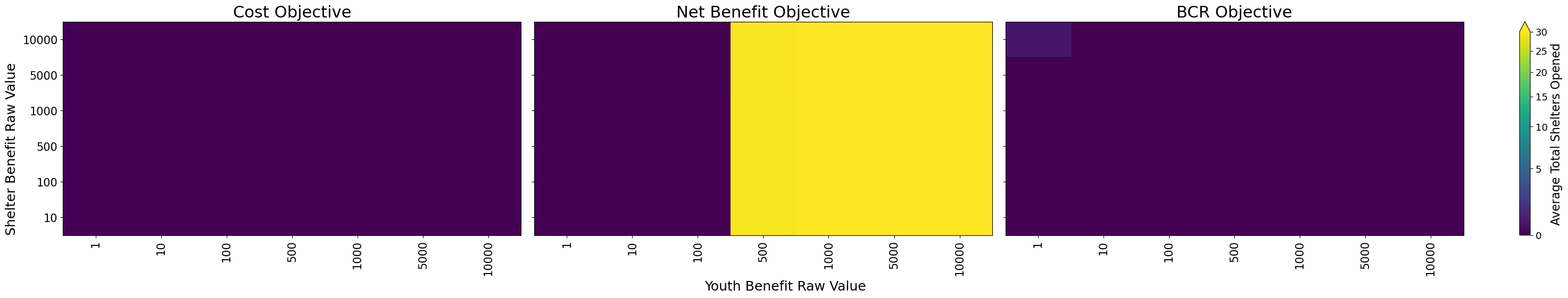}
    \caption{Comparison of conventional objectives with the BCR objective for $|Y| = 168$ youth and $\lambda = 2$. The axes report raw coefficient values for youth service benefits and shelter opening benefits.}\label{fig:NumericResults}
\end{figure}

\begin{figure}[H]
    \centering
    \includegraphics[width=\linewidth]{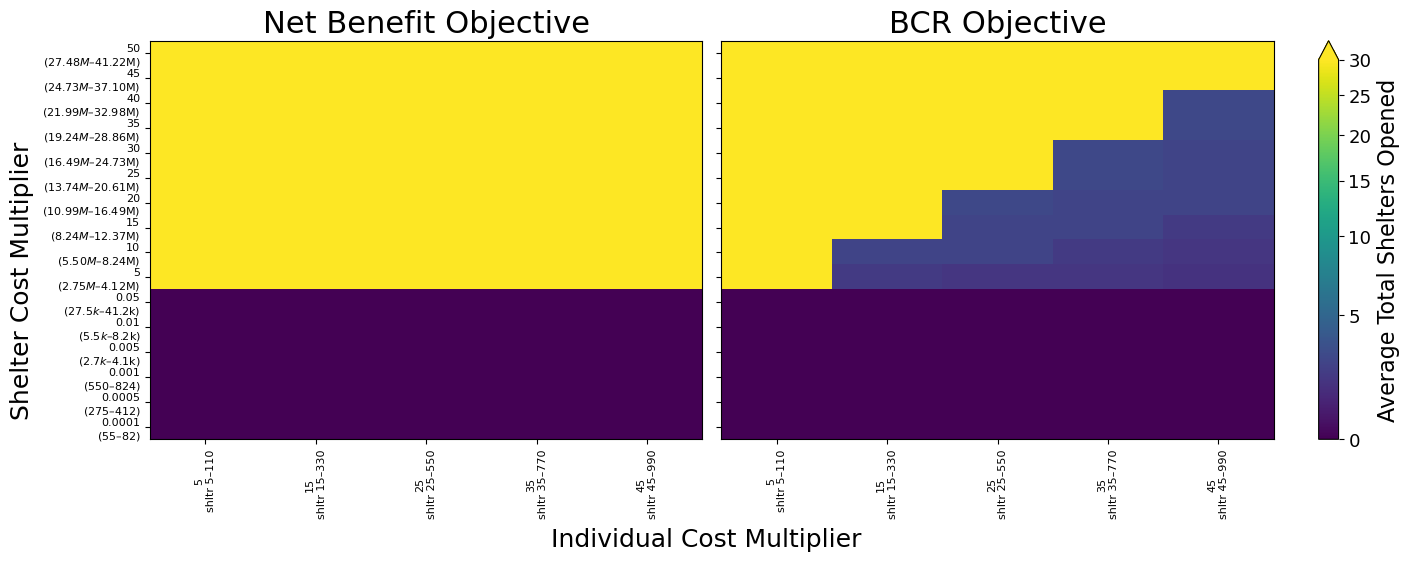}
    \caption{Comparison of conventional objectives with the BCR objective for $|Y| = 168$ youth and $\lambda = 2$. The axes report benefit multipliers applied to assignment costs and shelter opening costs. A more detailed view is provided in Appendix \ref{appendix:detailed_figs}. As the cost objective does not include benefits, this plot is omitted.}
    \label{fig:NumericResults_multiplier}
\end{figure}

This experiment compares the BCR objective with the two conventional objectives of cost minimization and net-benefit maximization, using $|Y| = 168$ youth and requiring at least $\lambda = 2$ shelter programs to open. % Here, we see that the ratio combines the advantages of benefits less cost maximization and cost minimization. 
Figures \ref{fig:NumericResults} and \ref{fig:NumericResults_multiplier} report the mean number of shelters opened across 10 sampled instances under two parameterizations of the benefit values. In Figure \ref{fig:NumericResults}, the horizontal axis represents raw youth service benefit values, $F_{y,s,i}$, while the vertical axis represents raw shelter opening benefit values, $\tilde{p}_{s}$. In contrast, Figure \ref{fig:NumericResults_multiplier} reports the same results under the multiplier parameterization, where the horizontal axis corresponds to the youth service benefit multiplier, $m$, and the vertical axis corresponds to the shelter opening benefit multiplier, $\tilde{m}$. 

Recall that, under the multiplier parameterization, $F_{y,s,i} = m \times r_{y,s,i}$ and $\tilde{p}_{s} = m' \times p_{s}$. Consequently, each point in Figure \ref{fig:NumericResults_multiplier} maps to a corresponding pair of raw benefit values in Figure  \ref{fig:NumericResults}, allowing the relationship between multiplier values and implied benefit magnitudes to be visualized. For example, a youth service multiplier of 15 multiplier of 15 produces raw values from 15 to 330, and $F_{y,s,i} = m \times r_{y,s,i}$ then this must imply that $r_{y,s,i}$  ranges from 1 to 22., while a shelter opening multiplier of 0.01 yields raw shelter opening benefits on the order of 27,500 to 41,200 due to the cost for opening a shelter $p_{s}$, ranging from 2.75 to 4.1 million.

As expected, Figure \ref{fig:NumericResults} shows that benefit values do not affect shelter opening decisions under cost minimization (the left panel), as the benefit function $B$ is not in the objective. By contrast, the net-benefit objective changes as the youth service benefit value increases (see the middle plot in Figure \ref{fig:NumericResults} and the left plot in Figure \ref{fig:NumericResults_multiplier}). We see that once the youth service benefit value becomes large enough for  benefits to outweigh costs, the net-benefit objective opens all shelters. On the right, we see the BCR objective does \textit{not} suffer from this, and the number of shelters opened remains at the minimum $\lambda = 2$ nearly always. %We note that the multipliers present in Figure \ref{fig:NumericResults_multiplier} can be mapped onto the raw values in Figure \ref{fig:NumericResults}, allowing one to see how the multiplier fits itself into the raw numeric values (for example, an youth multiplier of 15 gives actual youth benefit values on the order of 15 to 330, and a shelter multiplier of 0.01 gives actual shelter benefit values on the order of thousands). 

For example, when the shelter opening benefit (y axis) is $1000$ and the youth service benefit (x axis) is $500$, the net-benefit objective opens an average of all 30 shelters for the net-benefit objective, whereas the cost objective and the BCR objective only open the two shelters required by the action constraint. A more detailed view of this figure is provided in Appendix \ref{appendix:detailed_figs}. There, we provide the costs associated to each of these optimizers.

Figure \ref{fig:NumericResults_multiplier} shows that a youth service multiplier (x axis) of $25$ yields benefit values of $25-550$. Similarly, a shelter opening multiplier (y-axis) value of $5$ corresponds to benefit values $\tilde{p}_{s}$ between $2.75M$ and $4.1M$, as the cost of a shelter varies by location. This multiplier is in line the order-of-one multipliers discussed by~\cite{sroiValue}. This highlights the need for the BCR objective, as realistic shelter multipliers, combined with youth multipliers above order-of-one, allow us to determine optimization guided marginal action.

Lastly, we note the cost objective achieves the lowest average cost of $\$1,169,064$, while the BCR objective remains cost-conscious, with an average cost of $\$1,184,466$, and still identifies expansion decisions required by the action constraint.

\begin{figure}[ht]
    \centering
    \includegraphics[width=\linewidth]{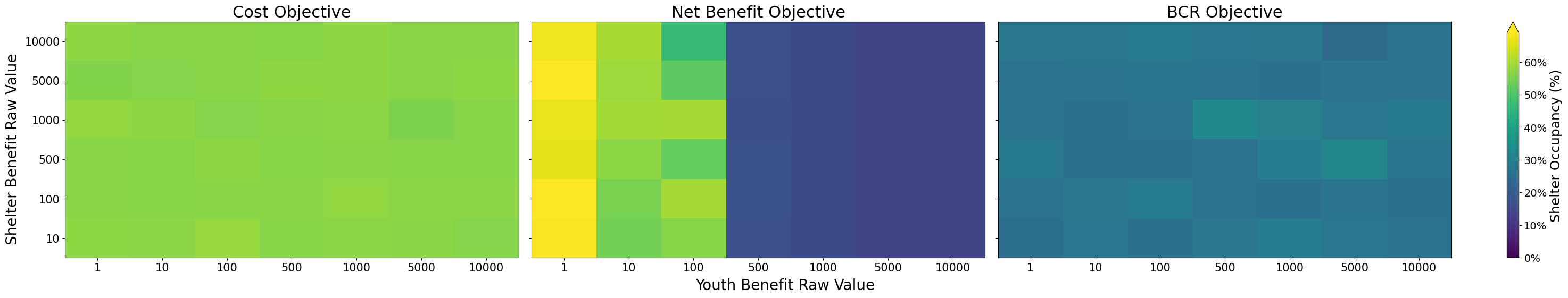}
     \caption{Average occupancy across operating shelter programs, including existing shelters and newly opened shelters, under the three objective functions for $|Y|=168$ youth and $\lambda=2$. A more detailed view is provided in Appendix \ref{appendix:detailed_figs}.}
    \label{fig:occupancy_changes}
\end{figure}

Compared with the net-benefit objective, the BCR objective changes the distribution of occupancy across the system while preserving a relatively low expansion cost. Figure \ref{fig:occupancy_changes} shows that, relative to the cost objective, the BCR objective increases new shelter occupancy from $23\%$ to $28\%$, a roughly $20\%$ relative increase, while increasing average cost by less than $2\%$. These results suggest that the cost objective relies more heavily on existing capacity, whereas the BCR objective makes greater use of newly opened shelter capacity. Although the cost and BCR objectives both open the two-shelter minimum and serve essentially the same number of youth in-house, they produce sharply different overall occupancy levels across all operating shelter capacity ($\approx 0.57$ vs.\ $\approx 0.27$). This occurs because the cost objective keeps existing shelter capacity densely used to minimize cost, whereas the BCR objective concentrates placements in higher benefit-per-dollar service slots, producing lower time-averaged occupancy even when total in-house service volume is similar. This lower average occupancy should be interpreted as preserved operational slack rather than higher utilization, providing available capacity if youth arrivals exceed expected arrivals.

% Overall, the BCR objective achieves nearly the same level cost while delivering a higher benefit–cost ratio, highlighting more efficient marginal expansion.

\subsubsection{\emph{Experiment 2: Sensitivity to Minimum Shelter Count}}\label{ss:sensitivityShelterLambda}

\looseness-1

There is a tradeoff between concentrating placements in fewer shelters and opening additional shelter programs to distribute service demand. Although placements can be concentrated in fewer shelter programs, opening additional shelters can reduce strain on existing capacity and increases system resilience. This tradeoff can help justify additional funding by quantifying the marginal benefit of opening new shelter programs. 

\begin{figure}[ht]
    \centering
    \includegraphics[width=0.5\linewidth]{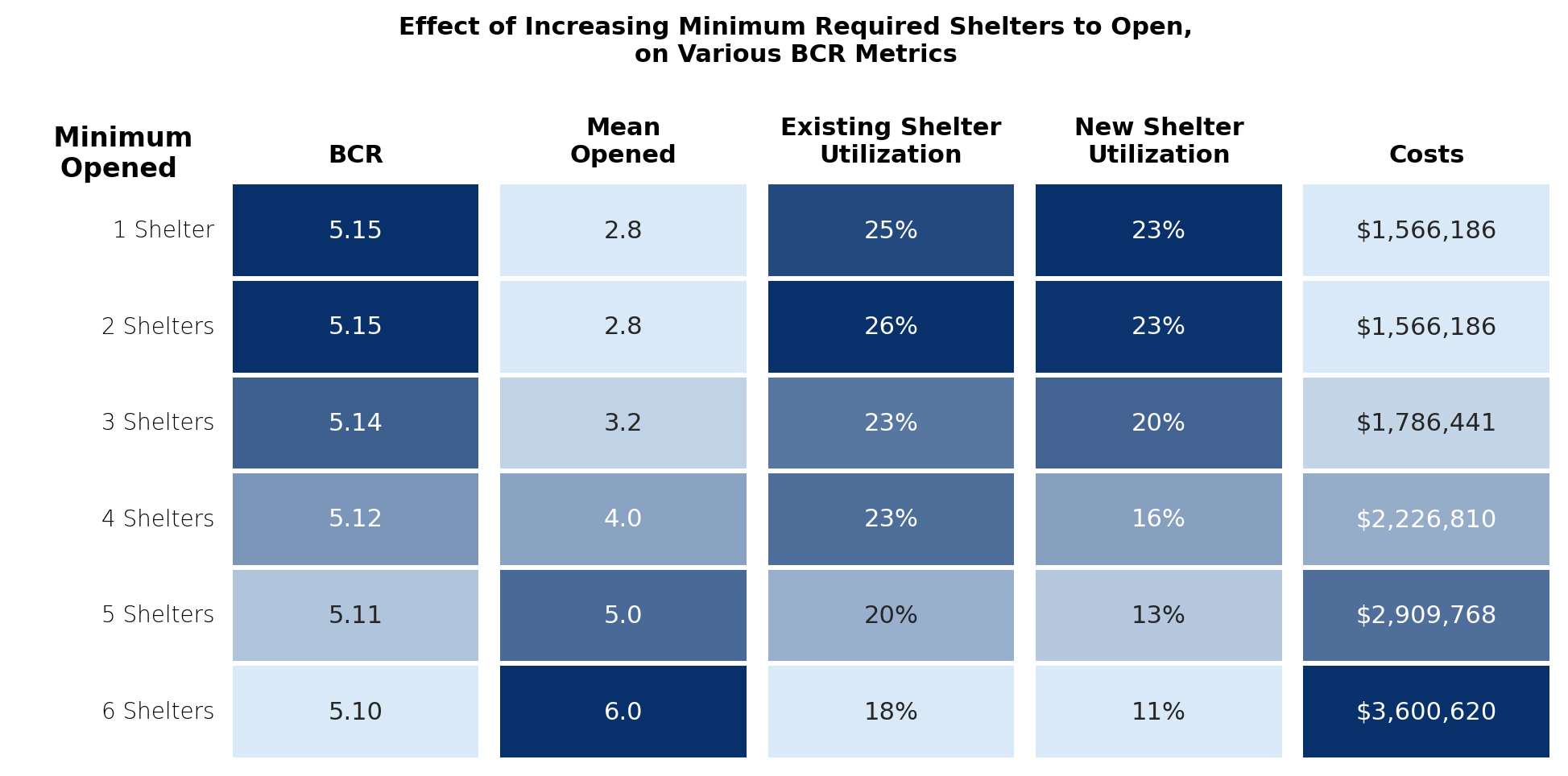}
    \caption{Effect of the minimum shelter-opening requirement $\lambda$ on BCR, number of shelters opened, occupancy, and cost under the BCR objective. The BCR objective is monotonically nonincreasing as $\lambda$ increases. }
    \label{fig:bcr_change}
\end{figure}

Holding data fixed, with a shelter multiplier of $5$ and a youth service multiplier of $15$, we optimize the BCR while varying the minimum number of shelter programs to open, $\lambda$.
Figure \ref{fig:bcr_change} shows how the BCR changes as $\lambda$ increases. This illustrates the role of action constraints in inducing movement away from the status quo.
Without action constraints, the BCR optimum may recommend preserving the status quo (``do nothing'') option; increasing $\lambda$ instead requires the model to identify controlled expansion plans that reduce strain on existing shelters and avoid overloading newly placed facilities.

\begin{figure}[ht]
    \centering
    \includegraphics[width=0.5\linewidth]{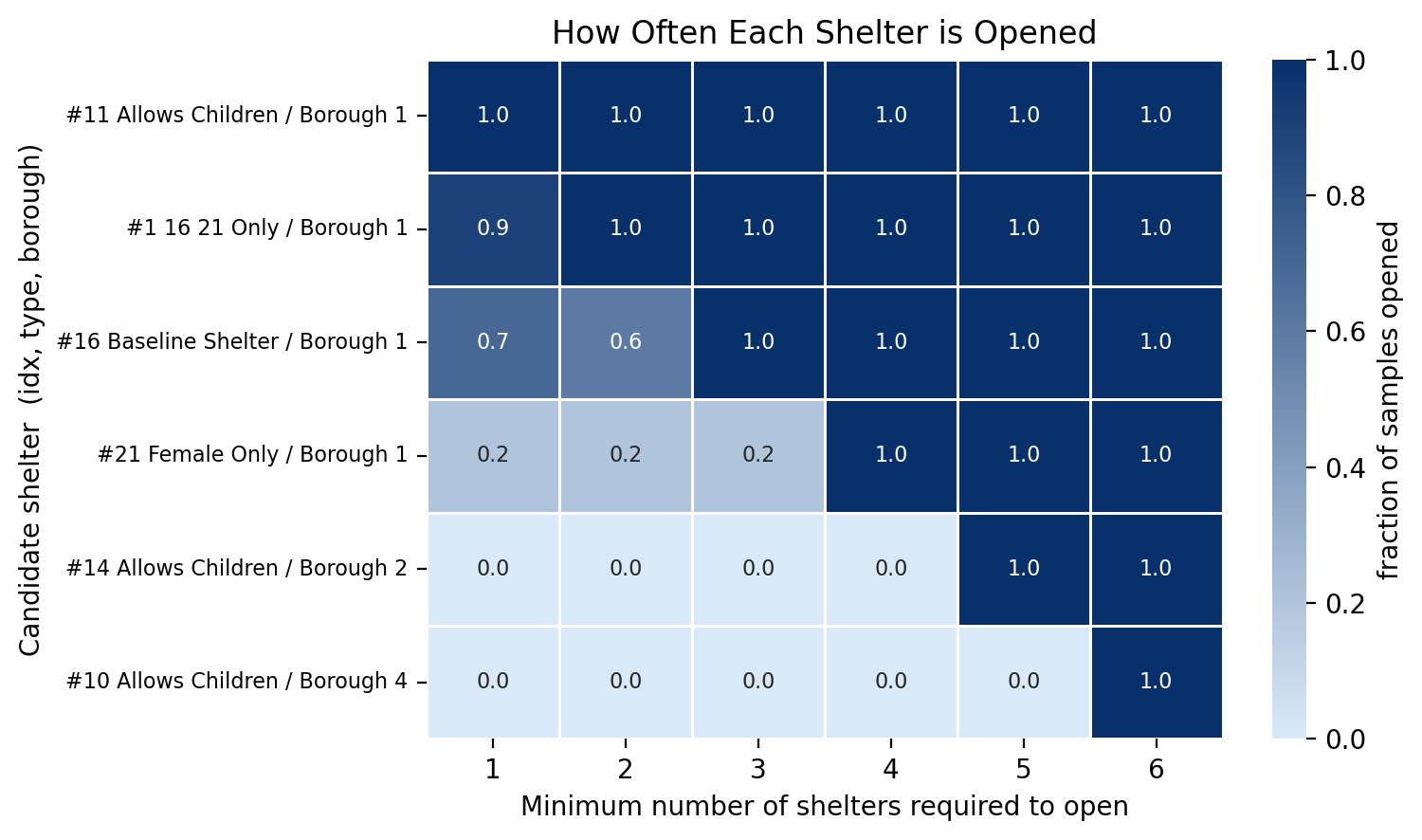}
    \caption{For each value $\lambda$, the fraction of samples (out of 10) for each shelter that are opened across runs under the BCR objective.}
    \label{fig:type_selection}
\end{figure}

Because existing shelters are modeled as operating at near capacity, we focus on the occupancy of considering the availability of the existing programs (that have limited capacities) as well as newly opened programs when evaluating how additional openings are used. 
Figure \ref{fig:type_selection} reports which candidate shelter types are opened across sampled instances for each value of $\lambda$.

As $\lambda$ increases, the BCR is monotone nonincreasing, while stable core of high-priority shelter programs opens first (e.g. at $\lambda=1$) and additional programs enter the solution at higher required levels of opening. The results show a consistent ordering of shelter types, with some types selected reliably at low values of $\lambda$ and others added only as the required number of openings increases.
This pattern suggests a phased deployment plan in which $\lambda$ is increased gradually to identify which new shelter programs should be added first, and which should be added later. This provides practitioners with a demand-planning approach for evaluating when additional capacity needs might justify increasing the required number of shelter openings.

Figure \ref{fig:type_selection} reports, for each value of $\lambda$, the fraction of the 10 sampled instances in which each shelter program type is opened under the BCR objective. These frequencies reveal a stable priority ordering: a core set of shelter types is selected first, and additional types added as the required number of openings increases. 

\subsubsection{\emph{Experiment 3: Adaptive Heuristic}} \label{sec:ScabailityModel}
% explain the results 
We consider three instance sizes: $|Y|=168$, $|Y| = 291$, and $|Y| = 717$, where the larger instances correspond to broader subpopulations of youth seeking additional housing service types; see Appendix \ref{dataDetailsAppendix} for details. For $|Y|=168$, the exact BCR model solves to the specified optimality tolerance; for larger instances, we apply the batch heuristic described in Section \ref{sec:core} with $\lambda =2$ over 10 samples and benchmark against the exact net-benefit solution computed on the full set of youth with a batch size of $56$ youth per iteration. We use the net-benefit objective $B - C$ because it solves quickly and represents a conventional decision-making criterion. Table~\ref{tab:heuristic} reports the resulting heuristic and benchmark performance.

\begin{figure}[H]
    \centering
    \includegraphics[width=1\linewidth]{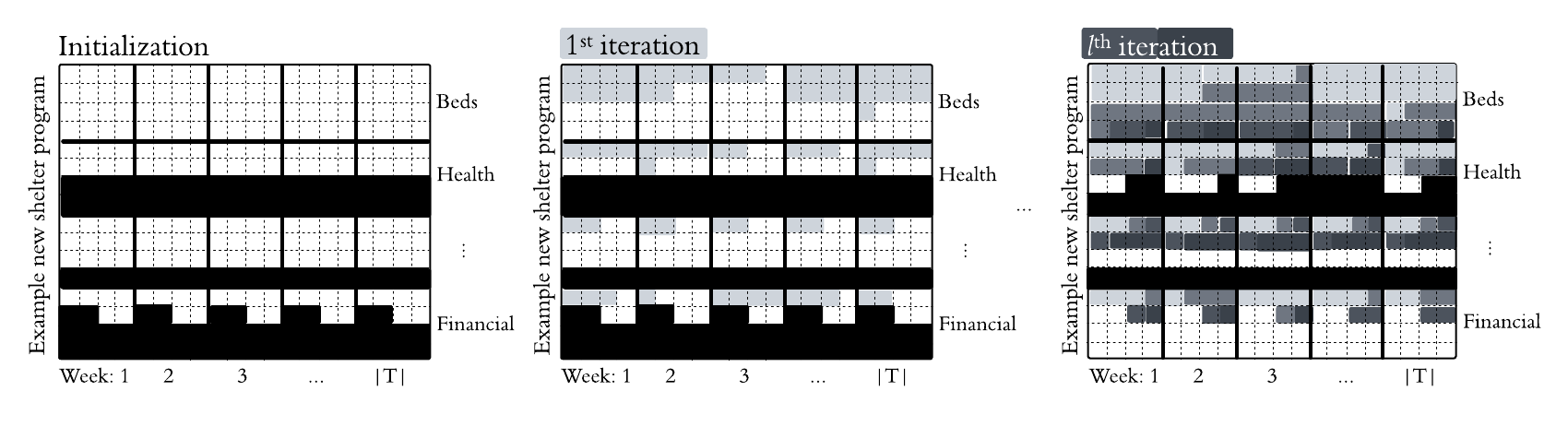}
    \caption{Illustration of successive CORE iterations. Each panel shows the assignments considered in a batch solve, with committed decisions carried forward as the heuristic updates the ensuing instance.}
    \label{fig:heuristic}
\end{figure}

Figure \ref{fig:heuristic} visualizes the CORE heuristic, in which each iteration solves a smaller batch model efficiently before committing selected decisions and proceeding to the next batch. The greyscale per cell indicate which slots are full, with black being completely occupied.

For $|Y| = 168$, the heuristic nearly matches the exact BCR optimum, achieving a BCR of 5.13 versus 5.15 and opening 3.1 versus 2.8 shelters on average across the 10 samples. For $|Y| = 291$, the heuristic opens 8.1 shelters compared with 30 under the exact net-benefit benchmark, approximately $73\%$ fewer, while achieving a nearly identical BCR, 5.09 versus 5.08. Both solutions result in zero  referrals.
For $|Y| = 717$, the heuristic opens 24.2 shelters compared with 30 under the exact net-benefit benchmark, approximately $19\%$ fewer, while maintaining a similar BCR, 5.13 versus 5.19.

  \begin{table}[H]
  \centering
  \renewcommand{\arraystretch}{0.7}
  \begin{tabular}{llcccc}
  \toprule
  {$|Y|$} & \textbf{Method} & \textbf{Shelters opened} & \textbf{BCR} & \textbf{Referrals} & \textbf{Runtime (s)} \\
  \midrule
  \multirow{2}{*}{168} & Exact     & \textbf{2.8}  & \textbf{5.15} & 0 & 13{,}319 \\
                       & Heuristic & 3.1           & 5.13          & 0 & \textbf{371} \\
  \cmidrule(lr){2-6}
  \multirow{2}{*}{291} & Exact     & 30.0          & 5.08          & 0 & \textbf{25} \\
                       & Heuristic & \textbf{8.1}  & \textbf{5.09} & 0 & 798 \\
  \cmidrule(lr){2-6}
  \multirow{2}{*}{717} & Exact     & 30.0          & \textbf{5.19} & 0 & \textbf{66} \\
                       & Heuristic & \textbf{24.2} & 5.13          & 0 & 1{,}829 \\
  \bottomrule
  \end{tabular}
  \caption{CORE heuristic performance averaged over 10 samples. The benchmark is the exact BCR solution for $|Y|=168$ and the net-benefit solution, for $|Y|=291$ and $|Y|=717$. Boldface indicates the preferred value within each instance size and metric.}
  \label{tab:heuristic}
  \end{table}
Overall, when exact optimization of the BCR objective is computationally prohibitive, the BCR-guided batches from Algorithm \ref{alg:core} maintain high service outcomes, with zero referrals in these experiments, while substantially limiting new program openings relative to the net-benefit benchmark, providing a practical and cost-sensitive deployment path.

\vspace{2mm}

\textbf{Key Insight: BCR Reveals Marginal Deployment of Resources} The BCR objective identifies high-value marginal resource-deployment decisions. Conventional objectives emphasize different extremes: cost minimization satisfies required action levels at low cost, while net-benefit maximization can open many shelters to reduce referrals, often with significantly higher expansion costs. 
In contrast, when action constraints require new shelter programs, the BCR objective weighs benefits relative to costs and prioritizes lower-cost expansions that reduce referrals while preserving cost discipline. 
In this sense, the BCR objective guides movement \textit{away} from the status quo by identifying the marginal actions that deliver the greatest improvement per unit cost under action constraints.

\section{Concluding Remarks}\label{keyTakeaways} % was Conclusions and Managerial Insights

\looseness-1
We develop a new algorithmic decision making framework that optimizes the benefit to cost (BCR) ratio to support public sector decision-making. While BCR comparisons are widely used, they are typically conducted manually and across only a limited set of alternatives. 
Our framework shows how a combinatorially large set of alternatives can be evaluated on a benefit to cost basis for optimization decision problems.
As advanced prescriptive analytics becomes increasingly available, this approach brings the power of data-driven optimization to BCR-based decision making in the public sector.

We demonstrate the framework through a compelling case study of critical capacity expansion decisions for runaway and homeless youth (RHY) shelter programs in New York City. 
% we introduce a decision context of expanding capacity to the homeless shelter system -- a public sector setting with extremely limited resources.
We estimate benefit and cost functions, introduce action constraints that induce movement away from the status quo, and optimize the resulting BCR objective using a variant of Newton's Method~\citep{Dinkelbach} that transforms the mixed-integer nonlinear optimization problem into a sequence of linearized mixed-integer problems.
Computational experiments show that the framework identifies high-leverage decisions for deploying marginal shelter program capacity.
The framework can be adapted to shelter systems in other communities or regions, provided appropriate data are available.

Optimizing the BCR objective provides finer granularity for identifying marginal improvements over the status quo than conventional net-benefit maximization or cost minimization.
Net-benefit maximization combines benefits and costs additively and can produce solutions that are relatively insensitive to implementation cost, thus impractical to adopt in resource constrained public sector contexts. Such approaches may also miss efficient marginal actions that deliver substantial improvement at relatively low cost. 
By contrast, the BCR formulation does not require benefits and costs to be combined as equally weighted additive terms.

\noindent \textbf{BCR Identifies Cost-Aware, High-Value Expansion Decisions.}
In our case study, BCR optimization contrasts with the net-benefit objective, which opens the maximum number of RHY shelters in many of the reported instances.
Conventional net-benefit objectives can obscure the marginal value of additional activity by emphasizing aggregate benefit net of aggregate cost. Our framework instead identifies ratio-maximizing near-term decisions and provides a foundation for iteratively developing expansion prescriptions over time.

The framework has important implications for public sector decision-making. Identifying  where to deploy marginal capacity can support strategic collaboration between researchers and decision-makers and generate effective, model-informed recommendations. This includes the opportunity to roll out marginal expansions sequentially, creating a path of improving feasible alternatives as the status quo evolves. Through action constraints, the decision maker can specify the minimal acceptable movement away from the status quo. In resource-constrained public sector contexts, this rollout approach offers an optimization-based strategy for phased capital allocation.

Several limitations suggest directions for future research.
To obtain tractable estimates of the benefits from satisfying action constraints, we assume that the corresponding partial returns remain constant.
Future research could incorporate richer dynamics among action variables and their partial returns, including the time required for shelters to open and become operational. Future work could also draw on principles of social return on investment, such as those in~\cite{sroi}, to estimate partial returns.
We also assume that youth have unrestricted access across the RHY shelter system; in practice, access may be geographically constrained by proximity to public transit lines.
Incorporating these access constraints could improve the operational realism of the model.
Similarly, future work could consider restrictions on which types of shelter programs can be opened in certain locations.

These limitations also point to many avenues for extension.
Analysts and decision-makers can iteratively evaluate alternative constraints to identify marginal actions that best support implementation, enabling co-development of data-driven models grounded in real-world decision contexts.
More broadly, these extensions underscore the potential of the proposed framework to deliver cost-effective, evidence-based decision support for public sector systems.
% For example, should decision-makers require shelters within a specific radius of a designated region $\mathcal{L}$, constraint~\eqref{eq:AtLeastlambdaOpenShelters} can be adapted to bound the number of  required locations within that region.

\section*{Data Availability Statement}
Data and code to repelicate studies is available upon reasonable request to the corresponding author.

\section*{Generative AI Statement}

Claude Opus 5.0 was utilized to create the figures in python for improved readability. Gen AI was not utilized in writing or editting the manuscript.

\if0\blind{
\section*{Acknowledgements}
We thank Dr. Meredith Dank and Andrea Hughes at New York University for their perspective on youth service provisions. We thank collaborators from the Coalition for Homeless Youth, Jamie Powlovich and Cole Giannone, and the New York City Mayor's Office under Eric Adams. 
This research was performed using computational resources supported by the Academic \& Research Computing group at Worcester Polytechnic Institute and the Research Computing Group at Northeastern University.
We are grateful to the National Science Foundation (Operations Engineering) grant CMMI-1935602 for their support as well as Grant DGE-2439018.} \fi

\bibliographystyle{chicago}
\spacingset{1}
\bibliography{bibliography}
\newpage
\section{Appendices}

\subsection{Model Tables}\label{mod:tables}
\begin{table}[H]
\centering

\renewcommand{\arraystretch}{0.7} % Decrease to compress rows
\begin{tabular}{ll}
\toprule
Set &  Definition \\ \midrule
$Y$ &  Set of youth \\
$S^{sq}$ & Set of existing organizations \\
$S^{sq-b}$ & Set of existing organizations that offer beds (and other services); ($S^{sq-b} \subseteq S^{sq}$) \\ 
$S^{sq-r}$ & Set of existing service referral organizations that do not offer beds; ($S^{sq-r} \subseteq S^{sq}$) \\
$S^{new}$ &  Set of potential new organizations  \\
$S^{service\ gap}$ &  Auxiliary set to identify youth unable to be served elsewhere \\
$S$ &  Set of organizations; ($S = S^{sq-b} \cup S^{sq-r} \cup S^{new} \cup S^{service\ gap}$) \\
$I$ &  Set of services offered by organizations \\
$D$ & Set of youth demographics\\
$T$ &  Set of discrete time periods\\ 
% $\mathcal{L}$ & $\ell$ & Set of locations where new shelter programs may exist \\
\bottomrule
\end{tabular}
\caption{Sets used in the marginal deployment optimization model.}\label{tab:modelSets}
\end{table}
\looseness-1

\begin{table}[H]
\centering

\renewcommand{\arraystretch}{0.7} % Decrease to compress rows
\begin{tabular}{lll}
\toprule 
Symbol & Domain & Definition \\ \midrule 
$\mathcal{K}_{s}$ & $\ZZ_{+}$ & Minimum number of youth needed to open $s \in S^{new}$\\ 
$\mu_{s,i}$    & $\ZZ_{+}$        & Maximum capacity of organization $s \in S^{sq} \cup S^{new}$ for service $i \in I$ \\ 
$c_{s,i}^t$   & $\ZZ_{+}$          & Capacity of organization $s \in S^{sq} \cup S^{new}$ for service $i \in I$ at time $t \in T$ \\
$\gamma_{s,i}$ & $\RR_{+}$ & Cost of expanding service $i \in I$ at organization $s \in S^{sq}$ by 1 unit    \\ 
$p_{s}$ &    $\RR_{+}$   & Cost of opening organization $s \in S^{new}$        \\ 
$r_{y,s,i}$& $\RR_{+}$       & Cost of assigning youth $y \in Y$ for service $i \in I$ at organization$s \in S$     \\
$\sigma_{s,i} $&$\{  0,1\}$             & $1$ If service $i \in I$ offered by organization $s \in S$ \\
%$\beta_{s}$ &$\{ 0,1 \}$   & Demographic profile (binary vector) of shelter program $s \in S$ \\
$\beta_{s,d}$ &$\{ 0,1 \}$   & 1 if organization $s\in S$ serves demographic profile $d \in D$ \\
$\tilde{p}_{s}$&$\RR_{+}$     & Benefit of opening organization $s \in S^{new}$  \\ 
$\rho $ &$\RR_{+}$& Benefit multiplier for $\tilde{c}(s)$ for $s \in S^{new}$ \\
$\lambda $ & $\ZZ_{+}$ & Minimum number of new organizations required to open \\
\bottomrule 
\end{tabular}
\caption{Organizational parameters used in the marginal deployment optimization model.}\label{tab:parameters2}
\end{table}

\begin{table}[H]
\centering

\renewcommand{\arraystretch}{0.7} % Decrease to compress rows
\begin{tabular}{lll}
\toprule
Symbol & Domain &Definition \\ \midrule
$\eta_{y,i}$ &$\{ 0,1 \}$   & $1$ If youth $y \in Y$ requests service $i \in I$ \\ 
%$\alpha_{y}$ &$\{ 0,1 \}$   & Demographic profile (binary vector) of youth $y \in Y$ \\ 
$\alpha_{y,d}$ &$\{ 0,1 \}$   & $1$ If youth $y \in Y$ has demographic $d\in D$ \\ 
$F_{y, s, i} $ & $ \RR_{+}$  & Youth benefits accrued (see Section \ref{ss:objective})  \\
$l_y  $ &$[T]$    & Arrival time of youth $y \in Y$ \\ 
$a_{y,i} $ &$[T]$  & Earliest start time of service $i \in I$ for youth $y \in Y$ \\ 
$b_{y,i}$ &$[T]$  & Latest start time of service $i \in I$ for youth $y \in Y$ \\ 
$d_{y,i}$  &$[T]$  & Duration of service $i \in I$ for youth $y \in Y$ \\ 
$f_{y,i} $   &$[T]$   & Frequency of service $i \in I$ for youth $y \in Y$ \\ 
\bottomrule
\end{tabular}
\caption{Youth parameters used in marginal deployment optimization model.}
\label{tab:parameters1}
\end{table}

\subsection{Variable Presolving and Cuts for Solvability}\label{app:presolve}
Due to the nature of the decision variables, our solver of choice is not able to handle large parameter spaces in multiple dimensions at once. We employ \textit{presolving} techniques, which eliminate certain variables by fixing them at $0$ before optimization begins. Additionally, we implement the following highly effective cutting planes. First, we provide the cuts: 
\allowdisplaybreaks
\begin{subequations}
    \begin{align}
        \sum_{y \in Y} X^{t}_{y,s,i} \leq \nu_{s} M, & \quad \forall s \in S^{new}, i \in I, t \in T \label{con:x_off_new_shelter_closed} \\ 
        \sum_{y \in Y} \pi_{y,s} \leq \nu_{s} M, & \quad \forall s \in S^{new} \\ 
        \sum_{y \in Y} \sum_{i \in I} U_{y,s,i} \leq \nu_{s} M , & \quad \forall s \in S^{new} \\ 
        \sum_{i \in I} U_{y,s,i} \leq \pi_{y,s} M, & \quad \forall y \in Y, \forall s \in S \\ 
        \sum_{t \in T} \sum_{i \in I} X^{t}_{y,s,i} \leq \pi_{y,s} M ,& \quad \forall y \in Y, \forall s \in S 
    \end{align}
\end{subequations}
Furthermore, we  eliminate  variables that would violate Constraint sets \eqref{con:demographic_match} and  \eqref{con:not_assigned_outside_range} by not initializing the corresponding $\pi_{y,s}$ and $U_{y,s,i}$ variables. This strategic variable reduction leads to a significant decrease in problem size, enabling computational tractability. 

\subsection{Youth Data and New Shelter Profiles Creation Explanation}\label{dataDetailsAppendix}

A survey by the authors and a group of social scientists [blinded for peer review] 
%\citep{dank2024understanding}
was conducted between December 2021 and June 2023 resulting in information on housing and service needs for 384 Runaway and Homeless Youth (RHY) in New York City (NYC). Using the RHY demographic information (Table \ref{tab:rhy_dem}) and service needs (Table \ref{tab:rhy_service}) reported by respondents, representative ``youth profiles" and ``service need profiles" were generated for the model. The total number of youth considered in the model experiments was informed by estimates from the January 2021 NYC point-in-time counts \citep{NYC_Youth_Count_2021}. Using these data we examine different scenarios for the number of RHY seeking shelter including variations for those seeking transitional independent living (TIL), short term housing, and emergency housing. We consider the latter two as RHY seeking these services often transition to seeking TIL within a relatively short time horizon. Further detail on the data generation process are found below.

\begin{table}[H]
\centering 
\begin{tabular}{ll} \toprule 
        Demographic & Data Type \\ \midrule 
        Gender & Cis-gender male, Cis-gender female, LGBTQ+,  Other/Unknown \\ 
        Sexual Orientation & Heterosexual,  LGBTQ+, Other/Unknown \\ 
        Has Children & Yes, No \\ 
        Legal Status & Legal, Illegal \\ 
        Human Trafficking Risk & Yes, No \\ 
        Age & Under 18, 18-21, Over 21 \\ \bottomrule
\end{tabular} 
\caption{RHY Demographics. Each item is a binary parameter.}
\label{tab:rhy_dem}
\end{table}

\begin{table}[H]
\centering 
\begin{tabular}{lp{11cm}} \toprule 
        Needs & Types of Services  \\ \midrule 
        Beds  \\  Health& (1) Mental Health, (2) Medical/Dental Health, (3) Substance and Alcohol Treatment \\ Personal Development&(4) Education, (5) Employment, (6) Practical Assistance, (7) Service Coordination, (8) Life Skills \\ Housing Support & (9) Long Term Supportive Housing, (10) Crisis Intervention \\ Legal & (11) Legal Assistance \\ Parenting &(12) Parental/Childcare Assistance \\ Financial & (13) Financial Support \\ \bottomrule
\end{tabular} 
\caption{RHY Service Needs. Each is a binary parameter.}
\label{tab:rhy_service}
\end{table}

The youth demographic and service need profile data was generated in the following way:
    \begin{enumerate}
        \item Youth Profiles: A probability density (PDF) function of unique "youth profiles" was generated using the observed combination of key demographic information across the RHY surveys (e.g., Age, Gender, Sexual Orientation, etc.). 
        \item Service Need Profiles: For each unique ``youth profile'', a set of ``service need profiles'' was generated using a PDF for each of the 13 types of services. This assumes each service need is dependent on the ``youth profile'' but independent from other service needs. We note that, given some services overlap and broadly cover similar needs, we simplify the original set of 13 services down to a set of 7 (see Table \ref{tab:rhy_service}) within the optimization model to improve solvability.
        \item Generated Youth Data: For each data point generated (youth), a random sample from the ``youth profiles'' PDF was selected and subsequently using the associated ``service need profile" PDF a value for each service need was assigned.

    \end{enumerate}

    The selection of number of youth to generate was determined in the following way: 
    According to the 2021 NYC Point-in-Time count, there were an estimated total of 2,054 unstably housed RHY in NYC on a given night \citep{NYC_Youth_Count_2021}. RHY experience different levels of need for housing and therefore only a proportion of these RHY would be requiring the form of housing (TIL) the presented model considers. Therefore, using the observed rate for specific housing needs identified in the survey ($n = 318$) (Question 200: If you were seeking shelter, what type of shelter would you want?) we examined three variations:  
    \begin{enumerate} 
    \item TIL seeking youth: 8.18\% of observed RHY were seeking shelter at a TIL, resulting in 168 youth considered under this scenario. This results in the most conservative number, representing only RHY whom are currently seeking a TIL specific bed.
    \item TIL and emergency or crisis seeking youth: 14.15\% of observed RHY were seeking shelter at either a TIL or an emergency or crisis shelter, resulting in 291 youth considered under this scenario. This results in a more robust number, representing both RHY whom are currently seeking a TIL specific bed, and individuals who are currently seeking an emergency or crisis bed, but may seek a TIL bed in the near future. 
    \item TIL and long term housing seeking youth: 34.9\% of observed RHY were seeking shelter either at a TIL or for long term housing opportunities, resulting in 717 youth considered under this scenario. This results in a larger number of youth given many youth desire long term housing. However, they are either currently be in a TIL or are not yet eligible for long term housing opportunities. 
    \end{enumerate}

The other data that is needed for this model is computed as follows. For the remaining youth data, we have that: 
\begin{enumerate}
    \item Arrival time is a uniform random integer in the time horizon $\{ 1,\dots ,T \}$.
    \item The duration for each service is a random integer chosen uniformly so that service requests are between $1$ and $3$ weeks. This then informs the total duration. Any time the duration goes longer than the current system horizon we model, we merely truncate and stop at the final time. 
\end{enumerate}

For the existing shelter information: 
\begin{enumerate}
    \item Costs for referral assignment and expansion come from a triangular distribution based upon the existing shelter programs. This is of the form 
    \begin{align}
        &r_{s,i} = TRIA(a_{s,i}, b_{s,i}, c_{s,i}), & \forall s \in R, \,& i \in I \\ 
        &\gamma_{s,i} = TRIA(d_{s,i}, e_{s,i}, f_{s,i}), & \forall s \in S^{sq},\, & i \in I
    \end{align}
    where $a_{s,i}, b_{s,i}, c_{s,i}$ are the low, mode, and maximum values for the referral assignment and $d_{s,i}, e_{s,i}, f_{s,i}$ are the low, mode, and maximum value for the expansion cost for existing programs. $\mu_{s,i}$ is computed in a similar manner. For data, see the reproducibility report. 
    \item Capacities in each shelter for nonbed services also come from a triangular distribution. Capacities for bed services are determined from our collaborators. 
\end{enumerate}

Then, the new shelter profiles are created in the following procedure:
\begin{enumerate}
    \item From the survey data and the existing shelters, we constructed $6$ shelter profiles for demographics our shelters could consider: the ``average'' of what exists, male only, female only, 16-21 only, 22-23 only, and one that explicitly allows RHY parenting youth. 
    \item Then we created $5$ copies of these shelter profiles and randomly assigned a borough to each one to determine the cost of the shelter. 
    \item The number of beds and other service capacities are the average size, which come from the existing shelters. 
    \item The remaining parameters are then averaged from the existing shelters. This way, our analysis focuses primarily on the needed demographics versus costs of opening new shelters.
\end{enumerate}
% \subsection{Proof of Satisfying Algorithm Requirements}\label{proofOfDinkelbachSatisfaction}
% \eb{
% For our problem, we must show that the linear relaxation of our model is compact, which is equivalent to the domain being closed and bounded. We have that the domain is closed since it is comprised of inequality constraints (non strict) and equality constraints. This means that our constraint set contains all of it's limit points, thus making it closed. To see what we have a bounded set, note that we have a finite number of new shelters we are considering, and all of our variables are bounded below by $0$ or bounded above by problem data. Therefore, we have a closed and bounded constraint set, which is therefore compact.  
% }

\subsection{More Detailed Figures}\label{appendix:detailed_figs}
\begin{sidewaysfigure}
    \centering
    \includegraphics[width=\linewidth]{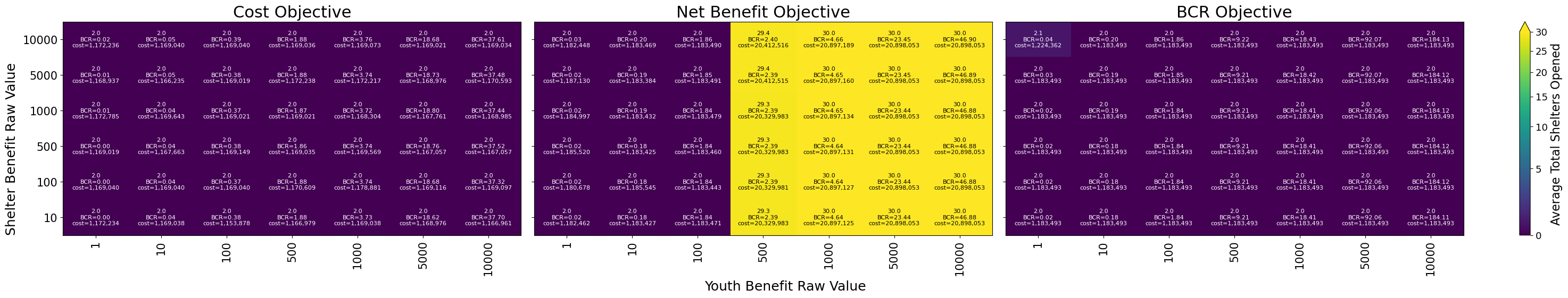}
    \caption{The same as \ref{fig:NumericResults} but we have visibility of the BCR, the average number of shelters opened, and the total costs.}
    \label{fig:placeholder}
\end{sidewaysfigure}

\begin{sidewaysfigure}
    \centering
    \includegraphics[width=\linewidth]{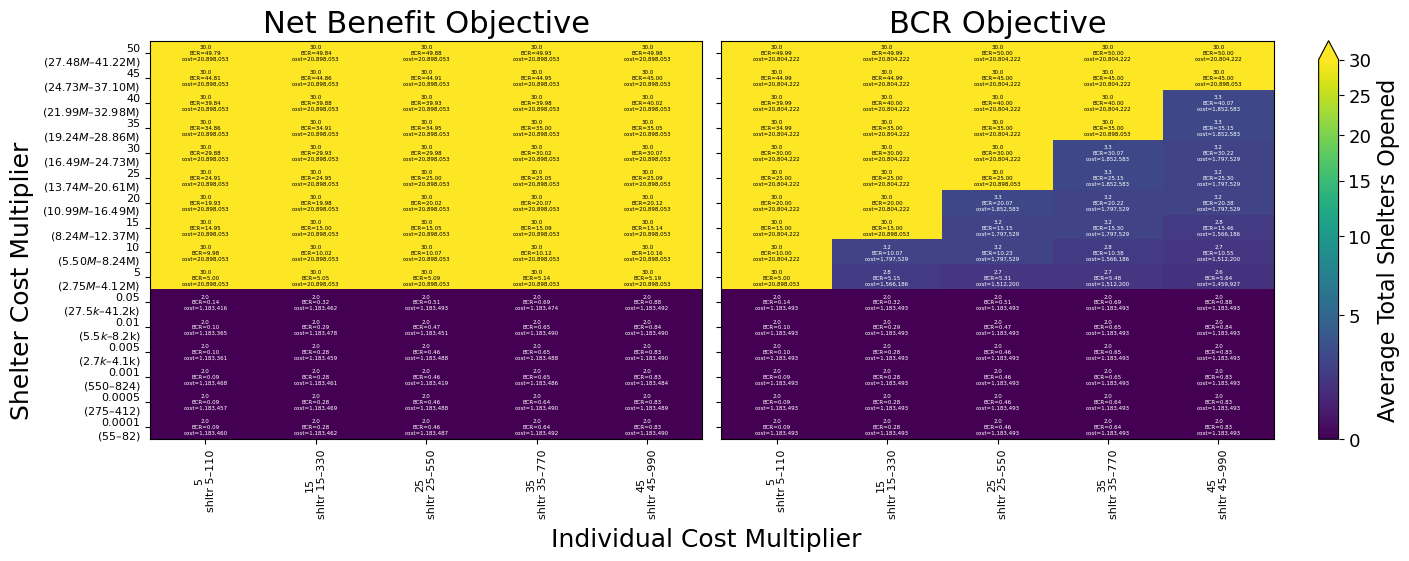}
    \caption{The same as \ref{fig:NumericResults_multiplier} but we have visibility of the BCR, the average number of shelters opened, and the total costs.}
    \label{fig:placeholder}
\end{sidewaysfigure}

\begin{sidewaysfigure}
    \centering
    \includegraphics[width=\linewidth]{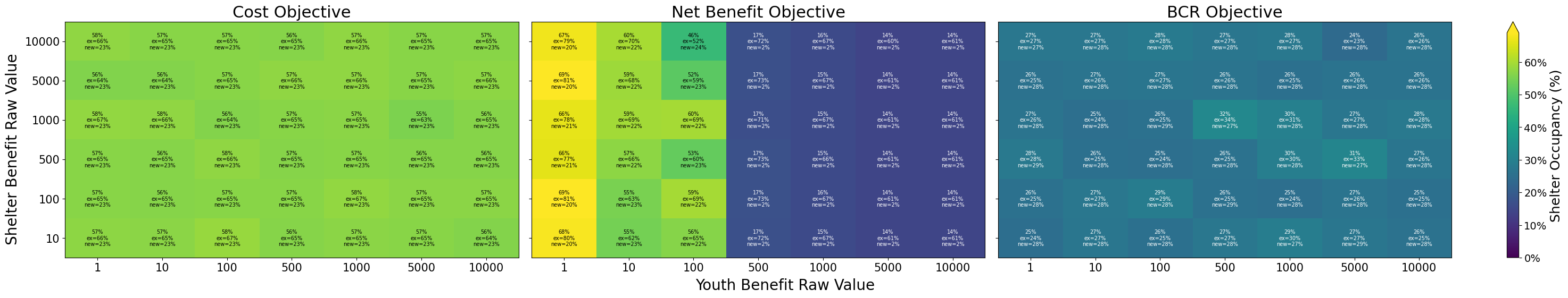}
    \caption{The same as \ref{fig:occupancy_changes} but we have visibility of the BCR, the average number of shelters opened, and the total costs.}
    \label{fig:placeholder}
\end{sidewaysfigure}

\end{document}